\documentclass{article}
\usepackage[utf8]{inputenc}
\usepackage{amsmath}
\usepackage{amsfonts}
\usepackage{amssymb}
\usepackage{graphicx,epstopdf,caption,subcaption}
\usepackage[left=2cm,right=2cm,top=2cm,bottom=2cm]{geometry}

\usepackage{algorithm,algorithmic,multirow}

\newenvironment{Proof}{\noindent{\sc Proof.}}{\qed}
\newtheorem{theorem}{Theorem}[section]

\newtheorem{definition}{Definition}[section]
\newtheorem{prop}{Proposition}[section]

\newcommand{\qed}{\hfill$\Box$\par\medskip}

\def\bhag#1{\noindent
\setcounter{equation}{0}
\section{#1}
}

\def\seqb{{\mathsf b}}

\def\RR{{\mathbb R}}
\def\ZZ{{\mathbb Z}}

\def\PPI{{{\rm I}\kern-1pt\Pi}}

\def\b #1;{{\bf #1}}
\def\x{{\bf x}}
\def\k{{\bf k}}

\def\m{{\bf m}}
\def\rr{\mathbf{r}^*}

\def\GG{\mathbb{G}}

\def\be{\begin{equation}}
\def\ee{\end{equation}}
\def\bea{\begin{eqnarray}}
\def\eea{\end{eqnarray}}
\def\eref#1{(\ref{#1})}
\def\disp{\displaystyle}
\def\cfn#1{\chi_{{}_{ #1}}}

\def\donchitre#1#2{\vskip 6.5cm\noindent
\parbox[t]{1in}{\special{eps:#1.eps x=6.5cm y=5.5cm}}
\hbox to 7cm{}\parbox[t]{0.0cm}{\special{eps:#2.eps x=6.5cm y=5.5cm}}}

\def\XX{{\mathbb X}}

\def\spans{\mbox{{\rm span }}}
\def\bs#1{{\boldsymbol{#1}}}
\def\GCN{{NHC}}
\def\LL{{MLL}}
\def\Percus{{MBO}}
\def\BC{{ES}}
\def\CC{{OS}}

\title {Representation of functions on big data associated with directed graphs}
\author{Charles K. Chui\thanks{Department of Statistics, Stanford  University, Stanford, CA 94305. The research of this author is supported by ARO Grant W911NF-15-1-0385.
\textsf{email:} ckchui@stanford.edu. }\ 
\  and
H.~N.~Mhaskar\thanks{Department of Mathematics, California Institute of Technology, Pasadena, CA 91125;
Institute of Mathematical Sciences, Claremont Graduate University, Claremont, CA 91711. The research of this author is supported in part by ARO Grant W911NF-15-1-0385.
\textsf{email:} hrushikesh.mhaskar@cgu.edu. } \
\ and
Xiaosheng Zhuang\thanks{Department of Mathematics, City University of Hong Kong, Tat Chee Avenue, Kowloon Tong,
Hong Kong. The research of this author is supported in part by the Research Grants Council of
Hong Kong (Project No. CityU 11304414) and City University of Hong Kong (Project No.: 7200462 and 7004445). \textsf{email:} xzhuang7@cityu.edu.hk. }
}
\date{}
\begin{document}
\maketitle
 
 \begin{abstract}
This paper is an extension of the previous work of Chui, Filbir, and Mhaskar (Appl. Comput. Harm. Anal.  38 (3) 2015:489-509), not only from numeric data to include non--numeric data as in that paper, but also from undirected graphs to directed graphs (called digraphs, for simplicity). Besides theoretical development, this paper introduces effective mathematical tools in terms of certain data--dependent orthogonal systems for function representation and analysis directly on the digraphs. In addition, this paper also includes algorithmic development and discussion of various experimental results on such data--sets as CORA, Proposition, and Wiki--votes.
 \end{abstract}

\bhag{Introduction}\label{intsect}

In this section, we first give a very brief summary on the recent progress of the manifold and (undiected) graph approaches for processing high-dimensional (numeric) data, and then discuss the need for directed graphs (called digraphs). 
We will also discuss the need for processing non-numeric data associated with digraphs, by using the endocrine network of the human body as an example. The objective of this paper is to develop a theory, along with a demonstration of some methods and algorithms, for the representation of functions on non-numeric data for the digraph paradigm, based on a  data-dependent orthogonal system, with associated filters, to be introduced in this paper. Since our approach is different from other studies in the literature, we will also give a toy example in this introduction section to illustrate the main idea of our approach. The organization of our presentation is outlined in the last paragraph of this section.

An earlier popular approach for processing high--dimensional numeric data is to consider that the data--set   lies near or on some (unknown) lower--dimensional manifold and to apply such mathematical tools as manifold learning, manifold Laplacian, diffusion maps and diffusion wavelets, to extract the data geometry and other data structures for data analysis. The interested reader is referred to the special issue \cite{achaspissue} of the journal, Applied and Computational Harmonic Analysis (ACHA), for some in--depth study in this research direction. 
In this regard, function approximation on such data--defined manifolds was also investigated in some depth (see, for example, \cite{mauropap, frankbern, modlpmz, eignet, heatkernframe}). 
On the other hand, since the discrete graph Laplacian well approximates the (continuous) manifold Laplacian (see \cite{singer} and the references therein), and since the subject of spectral graph theory (see \cite{chung1997spectral}) has already been a well established research area, it was almost immediately clear to at least a handful of researchers that perhaps high-dimensional data could be understood, analyzed, and processed  more fruitfully by associating the data directly with graphs, without embedding them into a lower--dimensional manifold. 
Indeed, many appealing aspects appearing in the analysis on the data--defined manifolds, such as the Hodge Laplacian,  various properties of the Green kernel, and useful inequalities, can be analyzed extensively in the context of spectral graph theory (see, for example, \cite{pang1993heat, smulders2004heat, chung2007heat, hodge_laplacian_lim2015}). 
In addition, function approximation on graphs has also been discussed in the literature. For instance, it is shown in \cite{friedman2004wave} that the solution of the wave equation corresponding to an edge--based Laplacian satisfies the property known as finite speed of wave propagation, and this, in turn, is equivalent to the Gaussian upper bound condition for small values of $t$ \cite{sikora2004riesz, frankbern}, so that the results in \cite{mauropap, tauberian} regarding function approximation and harmonic analysis are directly applicable. It is noted, however, that the associated graphs mentioned above are undirected graphs.

For big data analysis and processing, the associated graphs of interest are often directed graphs (or digraphs). For instance, digraphs are used effectively to model social networks, technological networks, biological and biomedical networks, as well as information networks \cite{newman2003structure}, with probably the most well-known example being the page--rank algorithm, where the nodes consist of urls of different web pages.
 Another example is the urban road network, where an intersection is represented by a vertex and a road section between adjacent intersections is denoted by an edge conforming to traffic flow direction on the road (e.g. \cite{han2012extended}). 
 The interested reader is referred to the special issue \cite{przulj2011introduction} of Internet Mathematics for an introduction to biological networks, with two papers \cite{elberfeld2011approximability, crofts2011googling} dealing directly with the question of finding the correct digraph structures. For biomedical informatics, the most important big data are arguably those of the human physiology; namely, the data generated by the physical, mechanical, and bio-chemical functions of humans and of human organs and systems. In particular, among the human body network systems, the two major ones, being the endocrine and nervous systems, both of which can be viewed as digraphs. While it is perhaps easier to understand, even by the layman, that our nervous system is a complex network of nerves and cells that carry messages from the brain and spinal cord to various parts of the body, the endocrine network is far more complicated. A very brief description of the endocrine system is that it is made up of a network of glands that secrete various chemical signals, called hormones, that travel in our blood vessels to regulate our bodily functions, helping in the control of our growth and development, internal balance of our entire body system, body energy levels (or metabolism), reproduction, as well as response to stress and injury. We will return to this discussion after briefly introducing the concept of non-numeric data.

An advantage of the (undirected or directed) graph approach over the manifold approach is that both numeric and non-numeric data can be (either separately or together) dealt with directly. Non-numeric data are also called qualitative data or categorical data, since they are used to indicate the quality of an object from the observed data (usually by using a bar chart) and showing the various categories in which the object belongs from the data (usually by using a pie chart). The usual techniques for working with numerical data cannot be used for non-numeric data. For example, the diffusion matrix commonly used for representing numerical data as a graph cannot be constructed in the same way for non-numeric data.  

Returning to the above discussion of the human body network systems, there are ten major glands that constitute the endocrine network, including the pituitary gland (also called the master gland, located at the base of the brain) that controls and regulates hormone secretion of the other glands.  Among these major glands, the pair of adrenal glands, with one sitting atop each kidney, are essential for human life. However, until today, there is still no reliable method for acquiring and analyzing the adrenal hormone data. Blood and urine tests are still commonly used, at least for preliminary screening to establish the case. The acquired information is non-numeric, in that adrenal insufficiency is (usually) determined by 7 observations, namely: sodium level, potassium level, blood pressure, glucose level, aldosterone level, cortisol level and ACTH level, in terms of only 3 qualitative marks: low, high, and normal, with categorial classification depending on personal genetics and medical history. If necessary, ultrasound or X-ray imaging of the abdomen to view the adrenal glands to establish primary adrenal insufficiency (called Addison's disease), and perhaps followed by CT scan to view the size and shape of the pituitary gland, if adrenal insufficiency could be secondary. When non-numeric data are associated with (undirected or directed) graphs, the data are represented as information in the nodes of the graph. The nature of this information is not critical to the analysis, but is used only to determine the edge weights of the graph. In the second paragraph, we have already discussed the topic of approximation of functions on numeric data associated with graphs. Here we mention that representation of functions on non--numeric data has also caught some attention recently, for instance, by Smale and his collaborators  \cite{smaleimmunology2012}, in introducing a mathematical foundation of molecular immunology, arguing that the study of peptide binding to some appropriate alleles can be thought of as a problem of approximating an unknown function on certain strings that represent the relevant alleles.

In the study of data associated with graphs, we note that digraphs are much less studied than undirected graphs in the literature. In fact, to the best of our knowledge, all current approaches to digrahps involve, in essence, the construction of an undirected graph that captures different features of the underlying digraph (see, for example, the recent surveys \cite{malliaros2013clustering} by Malliaros and Vazirgiannis or \cite{jia2014latest} by Jia et. al.).  For example, the Hodge Laplacian of a digraph is a symmetric matrix \cite{hodge_laplacian_lim2015}, and the weighted adjacency matrix of the graph Laplacian introduced by Chung \cite{chung_directed_laplacian} for a digraph is given by a symmetric  matrix as well, although an asymmetric version, called dilaplacian, has been discussed by Li and Zhang in \cite{li2012digraph}. The concept of generalized Cheeger constant that plays an important role in graph partitioning algorithms, as introduced in \cite{li2012digraph} also utilizes a symmetrized version of the dilaplacian. In the current paper, we propose an alternative way to develop harmonic analysis on digraphs by extending the undireted graph approach from our earlier paper \cite{treepap}. Our main idea is to represent a digraph by an asymmetric matrix $W$ (such as the weighted adjacency matrix, the dilaplacian, etc.), and observe that if the singular values of $W$ are all distinct, then $W$ can be recuperated uniquely from the symmetric matrices (equivalently, weighted undirected graphs) $WW^*$ and $W^*W$, or by their degree reduced forms. When the singular values are not distinct, the matrix is still a limit of matrices with distinct singular values. In other words, our viewpoint is that a digraph is a pair of undirected graphs. In this way, we can apply the well-known techniques for analysis on (undirected) graphs for developing an analysis on digraphs -- in principle, showing in fact that analysis on digraphs is trivial, once it is developed for (undirected) graphs. Therefore, in this paper we will represent a digraph as two (undirected) graphs and apply the theory and methods developed in our paper \cite{treepap}.

\begin{figure}[h]
\begin{center}
\begin{minipage}{0.3\textwidth}
\includegraphics[width=\textwidth,height=0.9\textwidth]{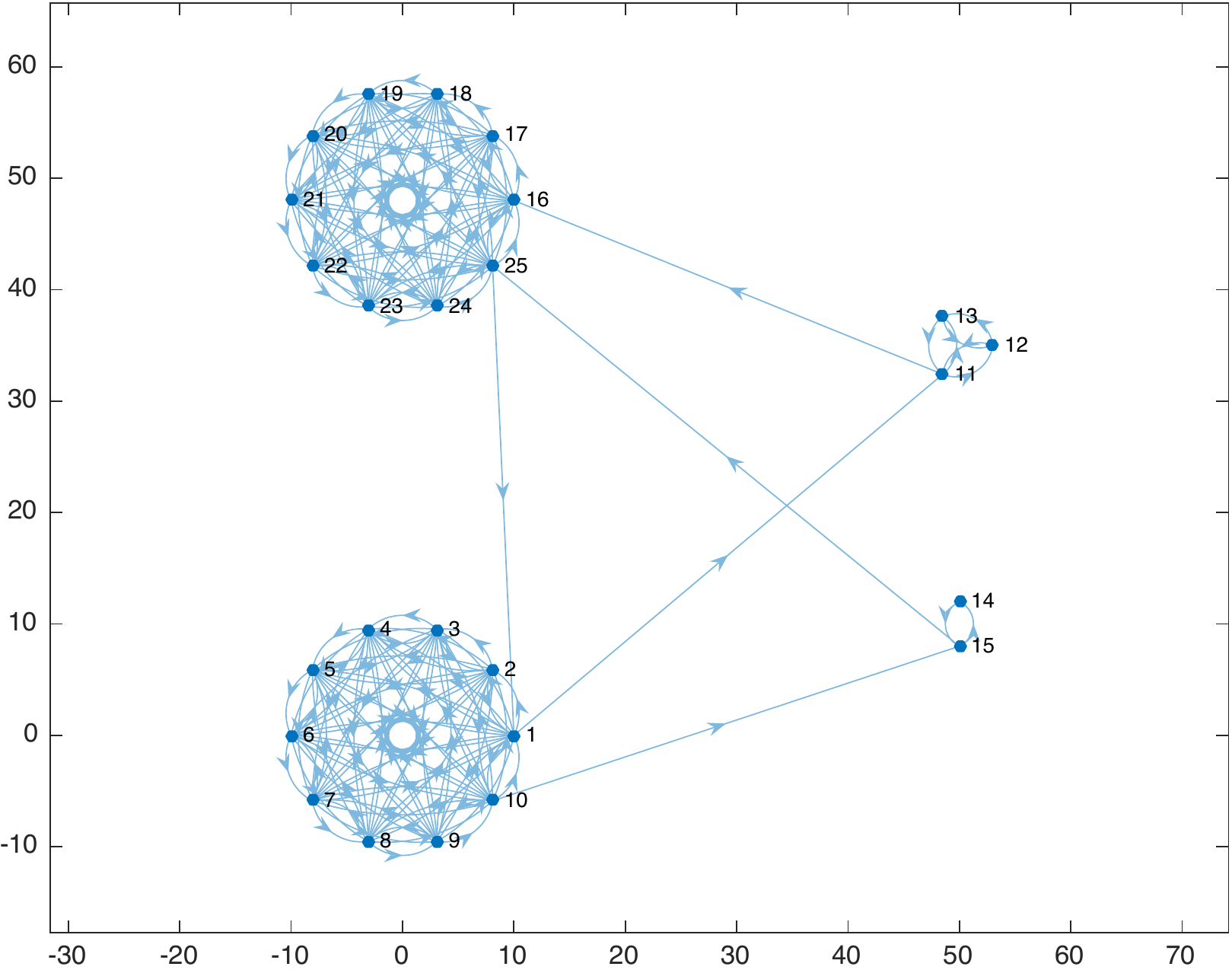} 
\end{minipage}
\begin{minipage}{0.3\textwidth}
\includegraphics[width=\textwidth,height=0.9\textwidth]{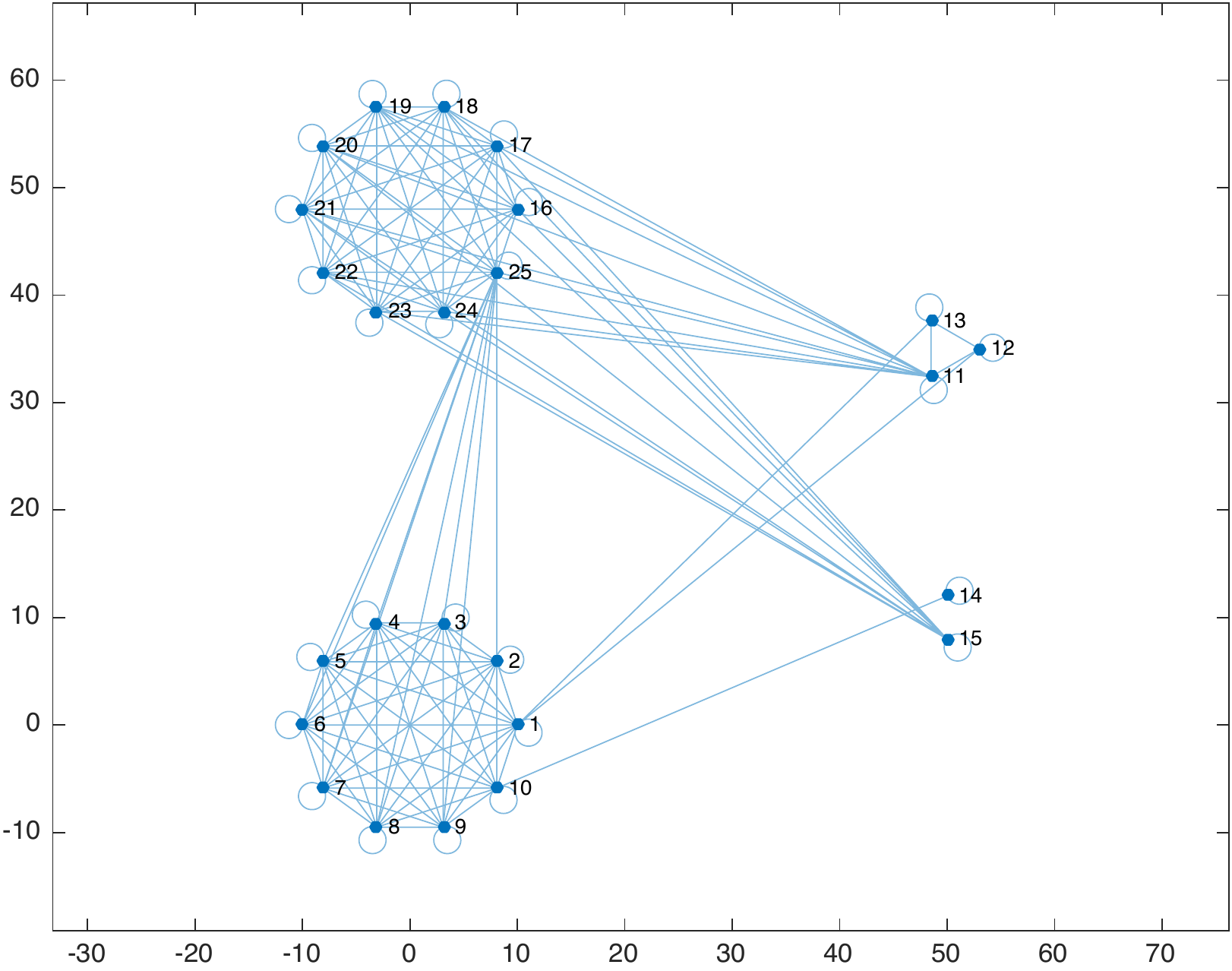} 
\end{minipage}
\begin{minipage}{0.3\textwidth}
\includegraphics[width=\textwidth,height=0.9\textwidth]{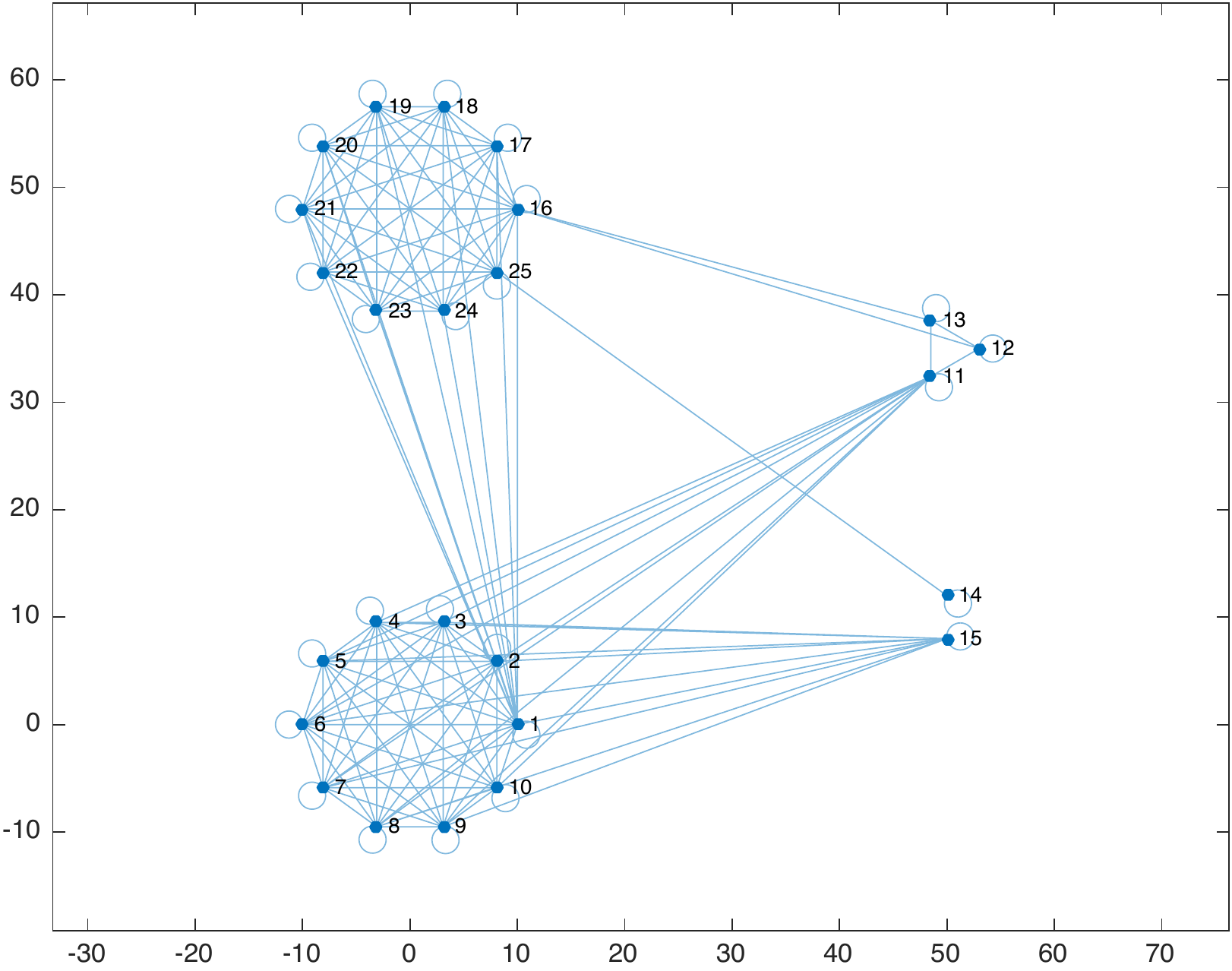} 
\end{minipage}

\end{center}
\caption{A simple, strongly connected digraph $W$ (left), and the equivalent pair of  graphs, $WW^*$ (middle) and $W^*W$ (right).}
\label{toygraphpict}
\end{figure}

\begin{figure}[h]
\begin{center}
\begin{minipage}{0.4\textwidth}
\includegraphics[width=\textwidth]{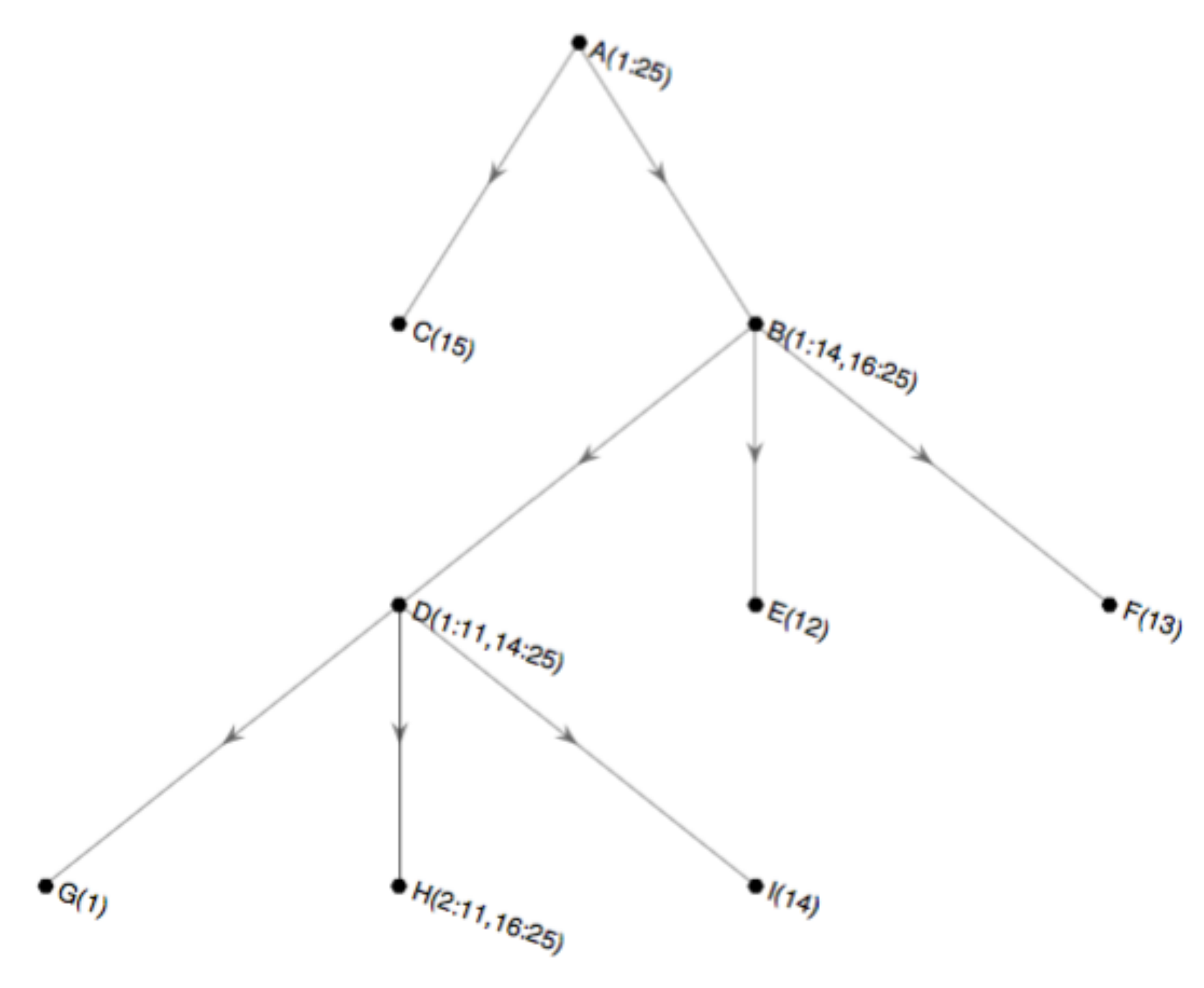} 
\end{minipage}
\begin{minipage}{0.4\textwidth}
\includegraphics[width=\textwidth]{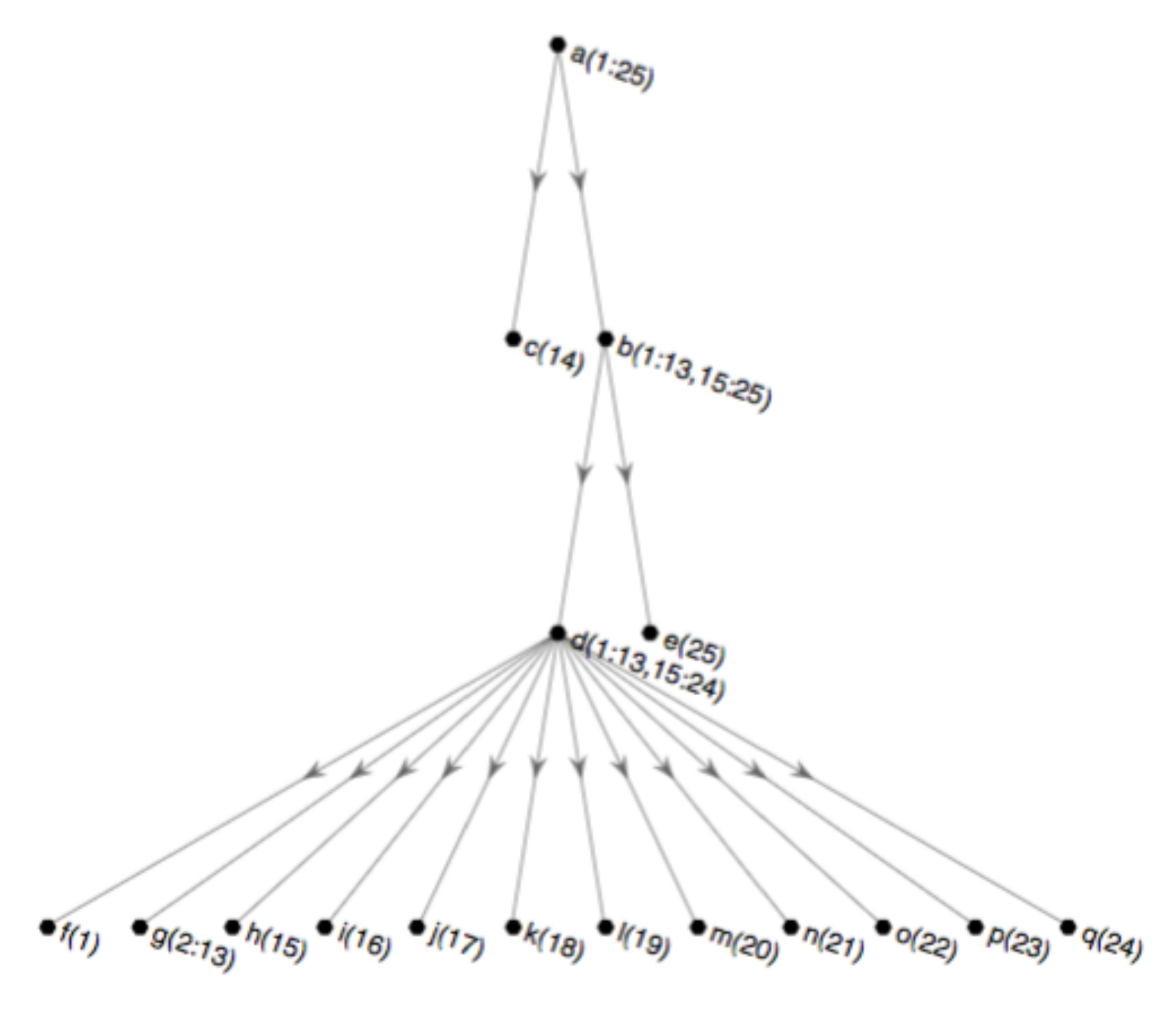} 
\end{minipage}
\begin{minipage}{0.4\textwidth}
 
\end{minipage}
\end{center}
\caption{The trees corresponding to the example in Figure~\ref{toygraphpict} using a variant of the  algorithm in \cite{dasgupta2010} with $WW^*$ on  left, and $W^*W$ on  right. Each node corresponds to either a cluster of nodes in $W$ as indicated in the parenthesis, or a node itself, also indicated in parenthesis. Thus, $q(24)$ is a leaf, corresponding to the node $24$ in $W$; $H(2:11,16:25)$ is a cluster with nodes $2$--$11$ and $16$--$25$. These nodes themselves are assumed to be the children of $H$. }
\label{toydaspict}
\end{figure}

\begin{figure}[h]
\begin{center}
\includegraphics[width=0.3\textwidth, height=0.27\textwidth]{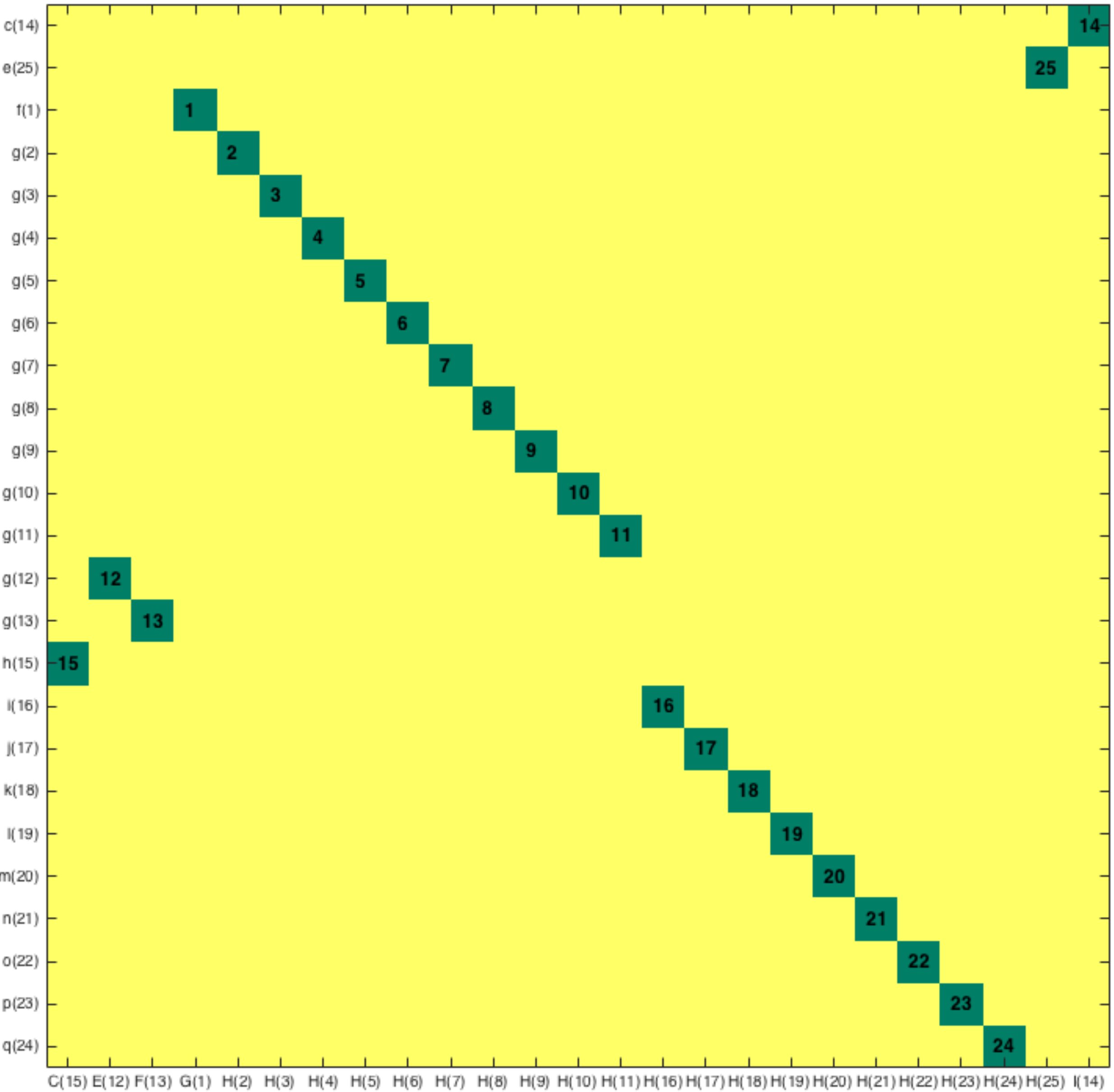}
\end{center}
\caption{The digraph $W$ of Figure~\ref{toygraphpict}  as rectangles in a partition of $I^2=[0,1)\times[0,1)$, according to the trees in Figure~\ref{toydaspict}.}
\label{toygraphmatpict}
\end{figure}

To demonstrate this idea, let us first consider a toy example with the digraph shown in Figure~\ref{toygraphpict}. Here, the weighted adjacency matrix $W$ is generated randomly but fixed throughout this example. In the general discussion of this paper, we will identify a digraph and its weighted adjacency matrix accordingly. If $W$ is the (weighted) adjacency matrix, we  apply a variant of the  algorithm for hierarchical clustering described by Chaudhury and 
Dasgupta \cite{dasgupta2010} for both $WW^*$ and $W^*W$, where the Euclidean distance is replaced by the graph distance on these graphs.  The resulting trees are shown in Figure~\ref{toydaspict}.  Although we do not show all the leaves of the two trees for the convenience of presentation, each of the trees corresponds to one node in the digraph $W$. Conversely, each node of $W$ appears as a leaf on each of the two trees. Using the edge weights of $W$, we can easily construct a filtration for each of the trees as described in \cite{treepap}, so that each of these leaves is a sub--interval of $[0,1)$. Suppose a node on $W$ appears as the interval $[a,b)$ on the tree corresponding to $WW^*$ and as the interval $[c,d)$ on the other tree. Then we will consider the node as the rectangle $[a,b)\times [c,d)\subseteq [0,1)\times[0,1)$, or according to convenience, any point on that rectangle. In particular, the digraph $W$ can now be viewed as a set of rectangles in a partition of $[0,1)\times[0,1)$, where each horizontal stripe as well as each vertical stripe contains at least one point of $W$. This is illustrated in Figure~\ref{toygraphmatpict}. In the sequel, $I^2$ will denote $[0,1)\times[0,1)$.

One major advantage of our approach is the following. It is observed in many examples, including the ones which we will study in this paper, that digraphs are typically highly disconnected. In contrast, the spectral theory for digraphs assumes a strongly connected (undirected) graph (e.g., \cite{chung_directed_laplacian, li2012digraph}). 
Our approach does not make any such assumption. We simply take the connected components of each of the two  (undirected) graphs as the children of the root for the tree corresponding to the graph, and use clustering on each of these connected components, enabling us to use spectral clustering if necessary. 
On the other hand, every point on the square $I^2$ does not correspond to a point on the digraph. Therefore, the theory of function approximation and harmonic analysis on digraphs in this paradigm must necessarily be totally data--driven, including data--dependent orthogonal systems and filters. In this paper, we will describe this theory in abstract. 

The outline of this paper is as follows. In Section~\ref{impsect}, we give examples of some of the algorithms used for hierarchical clustering
in order to represent a digraph as a pair of filtrations. We will also discuss standard criteria to evaluate the quality of this clustering.  These algorithms are tested in the case of three data--sets, an unlabelled data--set (Wiki-votes), a labeled data--set that is not hierarchically labeled (Proposition data), and a hierarchically labeled data--set (CORA). Each of these data--sets is non--numeric, and we make no effort to create a numeric data set out of them. The results are reported  in Section~\ref{datasetsect}. It is not our intention to investigate either the data--sets or the algorithms in their own right, but only to demonstrate that the choice of the algorithm can lead to a different structure for the digraph with vertices represented as elements of $I^2$. Therefore, unlike the study in classical harmonic analysis, not only the orthogonal system on the digraph, but also the very notion of smoothness and the various filters, must necessarily depend upon the data as well as the particular structure of the digraph obtained via the clustering algorithms. The theory of function approximation and harmonic analysis will be developed in Section~\ref{harmanalsect}.



\bhag{Implementation and testing}\label{impsect}

Our paper \cite{treepap} is motivated in part by the observation that a string as in \cite{smaleimmunology2012} can be represented via arithmetic coding as a vertex on a tree. Our work is motivated also by  the work of Coifman, Gavish, and Nadler \cite{gavish2012sampling, gavish2010multiscale}. Their approach starts with clustering the vertices of the graph into several ``folders'', followed by organizing these folders as another weighted graph, and repeating this process till only one folder remains. This organization generates a tree of sub-folders, sub-sub folders, and so on, till the lowest level that consists of the cluster-folders of the original vertices. In the paper \cite{treepap}, we have therefore assumed that the (undirected) graph has been converted to a tree using an appropriate clustering algorithm. In this section, we wish to extend this idea to digraphs. As explained in the introduction, a digraph can be viewed as a pair of undirected graphs. Clustering algorithms applied to each of these yield two corresponding trees, as well as a meaningful clustering of the digraph itself. The purpose of this section is to illustrate this concept using some concrete clustering algorithms and data sets.

After reviewing certain graph theory preliminaries, we review the algorithms we used for clustering (Sub--section~\ref{algsubsect}), as well as assessment tools for the quality of clustering (Sub--section~\ref{qualitysubsect}).

\subsection{Algorithms}\label{algsubsect}
For the convenience of the reader, we first review certain preliminaries about digraphs relevant to this paper.

A \textit{digraph} is an ordered pair  $(V, W)$, where $V$ is a  nonempty set, and $W :V\times V\to [0,\infty)$. Each element of $V$ is called a \textit{vertex} or a \textit{node}. 
If $u, v\in V$, then there is an \textit{edge from $u$ to $v$} with weight $W(u,v)$ if $W(u,v)>0$. This fact is often denoted by $(u,v)\in W$.
 The digraph is \textit{undirected} if $W$ is symmetric. 
The term graph (or undirected graph) refers to an undirected digraph.
A digraph $(V, W)$ is a \textit{tree} if there is a distinguished vertex $u^*$ (\textit{the root}) such that there is no edge from any other vertex to $u^*$, and for every $v\in V\setminus \{u^*\}$, there is a unique $u$ such that $(u,v)\in W$. 
The vertex $u$ is then called the \textit{parent} of $v$, and $v$ the \textit{child} of $u$. 
 We will follow the custom in computer science to treat the vertices as pointers to a record of information; for example, an entire file could be considered as the information stored in a vertex of some graph. For the convenience of exposition, we will often describe the vertex by the information it points to.

If $V$ is a finite set, then $W$ is represented by a matrix, called the \textit{weighted adjacency matrix}. 
If the values of $W$ are all in $\{0,1\}$ then $W$ is called an \textit{adjacency matrix}. In this section, we will assume $V$ to be finite, and denote the transpose of $W$ by $W^*$.

For a digraph $(V, W)$, the \textit{underlying undirected graph} is given by $(V, W_0)$, where $W_0=(W+W^*)/2$. If $I$ is the identity matrix of the same size as $W$, the \textit{extended graph} $(V, W_e)$ with  $W_e= I+W$ is the same graph as $(V,W)$ except for a new (or enhanced) self--loop inserted at each vertex. The \textit{pre--symmetrized (ES) graph} (respectively, \textit{post--symmetrized (OS) graph}) for $(V,W)$ is defined by $(V,W_{\BC})$ (respectively, $(V, W_{\CC}))$, where $W_{\BC}=W_eW^*_e$ and  $W_{\CC}=W^*_eW_e$. In the context of citation graphs, these have been called bibliographic coupling and co--citation graphs respectively \cite{diagraphcluster_ohiostate_2011}. 

If  $u, v\in V$, then a \textit{path from $u$ to $v$} is an ordered set $u_0=u, u_1,\cdots, u_n, u_{n+1}=v$ such that there is an edge between $u_i$ and $u_{i+1}$ for $i=0,\cdots, n$; the weight of this path is $\sum_{j=0}^n W(u_i, u_{i+1})$. The \textit{distance} from $u$ to $v$ is the minimum of the weights of all the paths from $u$ to $v$, defined to be $\infty$ if no such path exists. The distance matrix $d_W$ is the matrix whose $(u,v)$ entry is the distance from $u$ to $v$.

A (undirected) graph is \textit{connected} if for any $u, v\in V$, there exists a path from $u$ to $v$.  A digraph  is \textit{weakly connected} if the underlying undirected graph is connected. A \textit{weak component} of a digraph is a subgraph whose vertices form a maximal weakly connected subset of the vertices of the original digraph.   It is not difficult to show that if $(V,W)$ is weakly connected, then the {\BC} and {\CC} graphs for $(V,W)$ are connected (undirected) graphs. 

Each of the algorithms we discuss below have the  format described in Algorithm~\ref{genoutlinefig}, which we will call Twin--tree construction algorithm (TWT).

\begin{algorithm}[h]
\caption{\textbf{TWT}: A general top level description of the algorithms in this paper.}
\begin{algorithmic}[1]
\item[{\rm a)}] Input  a digraph $(V, W)$.
\item[{\rm b)}] Let $\{(V_j, W_j)\}_{j=1}^M$ be the weakly connected components of the extended graph $(V,W_e)$. With the digraph itself as the root, we construct two trees,  with the {\BC} (respectively, {\CC}) graphs for $(V_j, W_j)$, $j=1,\cdots, M$ as leaves. 
\item[{\rm c)}] Taking each of the leaves above as roots, we construct subtrees by applying various hierarchical clustering algorithms with the connected graphs represented by these leaves.

\item[{\rm d)}] The resulting trees are denoted by $\mathcal{T}_{{\BC}}^*$, respectively, $\mathcal{T}_{{\CC}}^*$. The symbol $\mathcal{T}^*$ will denote either of these.
\end{algorithmic}
\label{genoutlinefig}
\end{algorithm}

 In the sequel, we will describe our algorithms for a connected (undirected) graph; e.g., the leaf of the tree obtained in Step b of the algorithm TWT corresponding to the largest connected component of the corresponding undirected graph. Rather than complicating our notations, we will abuse the notation, and write $G=(V,W)$ for this graph, keeping in mind that in practice, this is really one of the leaves generated at Step~b of the algorithm in Figure~\ref{genoutlinefig}.  We denote the tree with root at $G$, resulting from the hierarchical clustering algorithm by $\mathcal{T}_G$.
 
Before describing the algorithms which we used extensively, we comment about some algorithms which we could not pursue vigorously.

The first algorithm to generate hierarchical tree structure from a connected graph is a variant of the  algorithm described by Chaudhury and Dasgupta in \cite{dasgupta2010}. This algorithm is developed primarily for clustering in high dimensional Euclidean spaces to achieve theoretically proven consistency results. In all the examples which we have studied in Section~\ref{datasetsect}, the data is not  numerical. Therefore, we replaced  the Euclidean distance  by the graph distance on $G$, as described in the introduction. However, we found this algorithm to be too slow for the examples.

The second one apparently highly cited algorithm is the Markov Clustering algorithm (MCL) \cite{van2001graph, shih2014component}. This is not
a hierarchical clustering algorithm, and therefore, we did not pursue this further. 

The other one is the MGL (Multiclass Ginzburg-Landau) algorithm described  in the paper \cite{arjuna2013}. Since it is similar to the MBO algorithm and MBO outperforms MGL in most of the cases, we therefore use only the MBO algorithm.

We now describe a set of three algorithms ({\GCN}, {\LL}, and {\Percus}) which we used in our examples as follows.

We developed an unsupervised, hierarchical clustering algorithm based on an idea described briefly  in \cite{gavish2012sampling, gavish2010multiscale}, that does not require the eigen-decomposition of a matrix. 
This variant is described in Algorithm~\ref{gcnalgorithmfig}. We will refer to this algorithm as {\GCN} (Non--spectral Hierarchical Clustering)  and
note that the algorithm is easy to modify for semi--supervised learning by choosing the initial centers to include the labeled data points.

\begin{algorithm}
\caption{\textbf{\GCN}: An unsupervised, hierarchical, eigen-decomposition free clustering algorithm.}
\begin{algorithmic}[1]
\item[{\rm a)}] {\bf Input}: undirected graph $G=(V,W)$. $K=(k_1,k_2,\ldots,k_L)$ with $1<k_1<k_2<\cdots<k_L<N$, where $N$ is the number of vertices in the graph.
\item[{\rm b)}] {\bf Output}: tree structure of level $0,1,\ldots, L+1$, where level 0 is the root $(V)$, level $L+1$ is the leaves of vertices, and in between are $k_l$ clusters at level $l$ for $l=1,\ldots,L$.
\item[{\rm c)}] {\bf Main Steps}:
\STATE Initialization: $\ell\leftarrow L$, $V_0\leftarrow V$, and $A_0\leftarrow W$.
\WHILE {$\ell>1$}
\STATE compute graph distance matrix $d_{A_0}$.
\WHILE{true}
\STATE randomly choose $k=k_\ell$ vertices $u_1,\ldots,u_k$ from $V_0$ as centers.
\STATE construct cluster $C_j$ for $j=1,\ldots,k$: $v\in V_0$ belongs to $C_j$ if $j = \mathrm{argmin}_{1\le i\le k}  d_{A_0}(u_i,v)$.
\STATE update the centers:  for each $C_j$, find a new center $u\in C_j$ such that $\sum_{v\in C_j} d_{A_0}(u,v)$ is minimal.
\STATE break if all centers remain the same.
\ENDWHILE
\STATE construct a new graph $G_1=(V_1,A_1)$ of $k$ vertices by the adjacent matrix $A_1$ of size $k\times k$ as follows: $A_1(i,j) = \sum_{i\in C_i, j\in C_j} A_0(i,j)$, $1\le i,j\le k$.
\STATE update $V_0\leftarrow V_1$, $A_0\leftarrow A_1$, and $\ell\leftarrow\ell-1$.
\ENDWHILE
\end{algorithmic}
\label{gcnalgorithmfig}
\end{algorithm}

At the other end of the spectrum, we used an  algorithm (Diffuse--Interphase Method)  described by Garcia--Cardona, et. al. in \cite{arjuna2013} which in turn is a modification of the well known MBO algorithm based on a graph Laplacian. This method can be used for hierarchical clustering only if the class labels are also organized hierarchically. Otherwise, we use this algorithm for the primary clustering, and use the coarse--graining ideas in \cite{lafonncut} to construct the remaining levels of the  tree bottom--up in an unsupervised manner. We will refer to this algorithm as {\Percus}.

In between the two, we used the algorithm described by Lafon and Lee in \cite{lafonncut}. This is also a  algorithm based on the graph Laplacian, but can be used both 
in the unsupervised setting (where the centers for clustering are chosen randomly) and in the semi--supervised setting (where the centers for clustering are chosen to be among  the training data). We will refer to this algorithm as {\LL} (Modified Lafon-Lee).

In both of the {\LL} and {\Percus}, we used the graph Laplacian. In our applications, it was not necessary to construct a diffusion matrix as in \cite{arjuna2013, lafonncut}. We only need the adjacency matrix of the graph as an input. In each case, the tree $\mathcal{T}_G$ has $4$ levels. The root is at level 0 containing all vertices and the leaves are vertices at level 3. We cluster all vertices to $k_2$ clusters at level 2 and   cluster them further into $k_1$ clusters at level 1 with $(k_1,k_2)$ preassigned.    Thus, a vertex $v$ of $\mathcal{T}_G$ at level $2$ is a cluster of the vertices in the graph $G$ which are children of $v$ in $\mathcal{T}_G$, and similarly, a vertex $u$ of $\mathcal{T}_G$ at level $1$ is a cluster comprising its children on the tree.

\subsection{Quality of clustering}\label{qualitysubsect}

It is clear that any harmonic analysis/function approximation scheme based on tree polynomials would depend upon the tree itself or equivalently on the quality of clustering used to generate the same. The objective of this paper is only to illustrate the concepts, not to point out an optimal clustering algorithm. Therefore, rather than using the usual measurement of accuracy of classification for evaluating our experiments, we will use measurements for the quality of clustering at different levels.

As explained in the introduction, each node on the digraph $W$ is interpreted as a rectangle contained in $[0,1)\times[0,1)$. The non--leaf nodes on the graphs would likewise be represented as rectangles as well, with each such node being the union of rectangles corresponding to its children. These non--leaf nodes at different levels will be considered as clusters at that level.

We will use two  measurements for the quality of clustering in digraphs in a hierarchical manner. To do so, we first make sure that the number of levels in the two trees corresponding to the digraph is the same. Suppose the tree $\mathcal{T}_{W_{\BC}}$ has $L$ levels and the tree $\mathcal{T}_{W_{\CC}}$ has $L'>L$ levels. Then we treat each  node at level $L$ in $\mathcal{T}_{W_{\BC}}$ as its own leftmost child, and continue this way until the tree  $\mathcal{T}_{W_{\BC}}$ has $L'$ levels as well. Equivalently,  since all the nodes  at level $L$ in $\mathcal{T}_{W_{\BC}}$ are leaves, any cluster at a level $>L$ is just a cluster according to $\mathcal{T}_{W_{\CC}}$. A cluster at level $\ell$ is then a rectangle in the partition of $I^2$ corresponding to the trees truncated at level $\ell$.

For unsupervised learning, we will use a measurement called modularity 
metric as described in \cite{malliaros2013clustering}. Various algorithms are recently designed to optimize this measure, for example, \cite{modularityopt2013}.
This metric is designed to measure the number of edges that lie within a cluster compared to the expected number of edges in a random digraph with the same in/out degree distribution. If $k_i^{\mbox{in}}$, $k_i^{\mbox{out}}$ represent the indegree, respectively the outdegree of node $i$ in $W$ (more precisely, the sum of weights on the incoming, respectively outgoing, edges at $i$), we assume that in a random digraph with the same connectivity, an edge from $i$ to $j$ will exist with probability $k_i^{\mbox{out}}k_j^{\mbox{in}}/m$, where $m$ is the total weight of the incoming/outgoing edges in the digraph; i.e., sum of the entries in $W$. Then the modularity metric introduced by Arenas, et. al. in \cite{arenas2007size} is defined by
\be\label{modularitydef}
\mathcal{M}=\frac{1}{m}\sum_{i,j}\left(W_{i,j}-\frac{k_i^{\mbox{out}}k_j^{\mbox{in}}}{m}\right)\delta(C_i,C_j),
\ee
where $\delta(C_i,C_j)$ is $1$ if the nodes $i$ and $j$ both belong to the same cluster $C=C_i=C_j$, and $0$ otherwise.
In our implementation of this metric, we will consider the nodes at each level of the trees as the clustering at that level.

  In the semi--supervised setting, we used the $F$--measure described in \cite{diagraphcluster_ohiostate_2011}. If $\{C_1,\cdots,C_M\}$ are the obtained clusters in the digraph from certain clustering algorithm, and $\{L_1,\cdots,L_n\}$ is a partition of the nodes according to the (ground-truth) class labels (i.e., $L_j$ is the set of all nodes in $W$ with the class label $j$), then one defines
$$
F(C_i)=2\max_{1\le j\le n}\frac{|C_i\cap L_j|}{|C_i|+|L_j|}.
$$
the (micro--averaged) $F$--measure is then defined by
\be\label{fmeasuredef}
\mathcal{F}=\frac{\sum_i |C_i|F(C_i)}{\sum_i |C_i|}.
\ee

In Sub-section~\ref{diharmsect}, we use the confusion matrix (see \eqref{def:cm}) to measure the approximation power of different algorithms using our framework.

\bhag{Data sets and results}\label{datasetsect}

We present our results for the (1) CORA data set  (2) Proposition data set, and the (3) Wiki--votes data set.

Each of these data sets contains only one large weakly connected component while others are of small size.  In our exposition, we  focus only on the largest weakly connected component as the leaf obtained in Step~b of the algorithm in Figure~\ref{genoutlinefig}.
 The same algorithms can be applied to the  other   weakly connected components. 
 If these components are too small, we may treat their vertices as the children of the $\mathcal{T}^*$--vertices corresponding to these components. 
 
By abuse of notation as before, let $G=(V,W)$ be the subgraph with respect to the largest weakly connected component, $W_{\BC}$, $W_{\CC}$ be the {\BC} (respectively, {\CC}) graphs of $(V, W)$, and $\mathcal{T}_{W_{\BC}}$ (respectively, $\mathcal{T}_{W_{\CC}}$) be the resulting $4$ level trees with $G$ at its root.

For each of the methods ({\GCN, \LL, \Percus}),  we randomly pick $p\%$ of the data as training data (semi-supervised learning (SSL) while 0\% means unsupervised learning (USL)) and perform the clustering algorithms.  For the method \Percus, we used $50$ eigenvectors, the time step is $0.01$, the stop criterion is $10^{-3}$, and the weight constant for the fidelity term is $50$.  For the algorithm \LL, we used $30$ significant eigenvectors and ``time parameter'' $t=1$.  We compute the modularity metric and $F$-measure for each of the levels as described in the introduction.  In view of the random choices of centers in both unsupervised (USL) and semi--supervised (SSL) settings, we computed these measurements for each given $p\%$ over 30 trials, and the modularity metric $\mathcal{M}$ and $F$-measure $\mathcal{F}$ are average over these 30 trials.  Note that for the Wiki-votes data set, we can only compute the modularity metric. 
  
\subsection{The data set CORA}\label{corasect}

We worked with the CORA research paper classification data set downloaded from https://people.cs.umass.edu/$\sim$mc-callum/data.html. The data set comprises a digraph with 225,026 publications as vertices, and edge from $i$ to $j$ means that paper $i$ cited paper $j$. These publications are from several areas of computer science and information theory. The subject area of each publication is given at two levels; e.g., 
artificial intelligence/vision, artificial intelligence/agent, artificial intelligence/DNP.  Out of the entire data set, only 28,135 are labeled. We considered only the subgraph  whose vertices are from this labeled data set. Altogether there are 70 classes at the most refined level, which are subgrouped into 10 classes, yielding a hierarchically labeled data set. There are 4,070 weakly connected components, and 22,985 strongly connected components, most of which are singletons. Thus, the digraph is highly disconnected.  The largest weakly connected component of $G$ contains 23,567 vertices while other weakly connected components contain at most 12 vertices. For each of $W_{\BC}$ and $W_{\CC}$ from the largest weakly connected component, we cluster all vertices to 70 clusters at level 2 and then further cluster them to 10 clusters  at level 1.

For this data set, we can perform both unsupervised learning and semi-supervised learning methods. The results are given in Table~\ref{corarestab:F-M}. 
From the table, in terms of the $F$-measure and modularity metric,  semi-supervised method {\Percus} performs better than the other two methods of {\GCN} and {\LL}, especially when the size of the training data is small ($\le30\%$).  
It is interesting to note that the best $F$--measure reported in \cite{diagraphcluster_ohiostate_2011} for this data set is $0.36$ at level $2$, while the the algorithm {\Percus} applied with 40\% training data yields a better $F$--measure both at levels $1$ and $2$.

\begin{table}[ht]
\begin{center}
\begin{tabular}{c|l|c|c|c|c|c|c|c|c|c|c}
\hline 
\hline
$\mathcal{F}$ & Trains (\%) & 0 (USL) &  10 & 20 & 30 & 40 & 50 & 60 & 70 & 80 & 90  \\
\hline
\multirow{2}{*}{\GCN} &
Level 2 (70)& 0.10& 0.13& 0.17& 0.23& 0.31& 0.41& 0.51& 0.62& 0.74& 0.87\\
& Level 1 (10)& 0.15& 0.49& 0.49& 0.54& 0.59& 0.65& 0.71& 0.77& 0.84& 0.92\\
\hline
\multirow{2}{*}{\LL} &
Level 2 (70)& {0.11}& 0.12& 0.15& 0.21& 0.29& 0.38& 0.49& 0.60& 0.73& 0.86\\
&Level 1 (10)& 0.42& 0.40& 0.40& 0.45& 0.48& 0.52& 0.58& 0.67& 0.77& 0.87\\
\hline
\multirow{2}{*}{\Percus} &
Level 2 (70)& N.A.& 0.21& 0.27& 0.33& 0.40& 0.49& 0.58& 0.68& 0.78& 0.89\\
&Level 1 (10)& N.A.& 0.54& 0.56& 0.61& 0.61& 0.66& 0.73& 0.79& 0.86& 0.93\\
\hline
\hline
$\mathcal{M}$ & Trains (\%) & 0 (USL) &  10 & 20 & 30 & 40 & 50 & 60 & 70 & 80 & 90  \\
\hline
\multirow{2}{*}{\GCN} &
Level 2 (70)& 0.34& 0.32& 0.32& 0.32& 0.34& 0.36& 0.40& 0.45& 0.49& 0.55\\
&Level 1 (10)& 0.35& 0.30& 0.33& 0.34& 0.36& 0.38& 0.42& 0.46& 0.50& 0.55\\
\hline
\multirow{2}{*}{LL} &
Level 2 (70)& 0.23& 0.23& 0.24& 0.26& 0.29& 0.32& 0.37& 0.42& 0.47& 0.55\\
&Level 1 (10)& 0.22& 0.24& 0.23& 0.25& 0.28& 0.31& 0.36& 0.40& 0.46& 0.53\\
\hline
\multirow{2}{*}{\Percus} &
Level 2 (70)& N.A. & 0.51& 0.48& 0.45& 0.44& 0.44& 0.45& 0.48& 0.51& 0.56\\
&Level 1 (10)& N.A. & 0.52& 0.49& 0.47& 0.48& 0.48& 0.50& 0.51& 0.53& 0.56\\
\hline
\hline
\end{tabular}
\end{center}
\caption{CORA Data: the $F$-measure $\mathcal{F}$ and modularity metric $\mathcal{M}$ using {\GCN, \LL},  and {\Percus} for given 0\% (USL method), 10\%, 20\%, \ldots, 90\% training data at level 2 ($k_2=70$ clusters) and level 1 ($k_1=10$ clusters), respectively. All results are average over 30 trials.}
\label{corarestab:F-M}
\end{table}

\subsection{The Proposition data set}\label{propdatasect}
This data set is described in detail in \cite{zhu2014tripartite, smith2013role}. The November 2012 California ballot contained 11 initiatives, or propositions, on a variety of issues, including state taxation, corrections, and food labelling among others. The  data consist of Twitter posts related to initiatives, grouped according to different propositions. For each proposition, the data is a directed graph with edge from $i$ to $j$ if the tweet originated from user $i$ to user $j$. The authors of \cite{smith2013role} have assigned an evaluation of the emotion of the sender with each tweet, called sentimental values, thereby creating a real valued label. 
At level 1, we  group the users into binary clusters by the sign of sentimental values (yes or no). At level 2, we  divide these further into 10 finer groups according to the strength of the sentimental values (strongly disagree to strongly agree). In such a way, we can construct a hierarchical labelling and all the algorithms can be used here as well. We choose the largest data set (Prop 37). Its largest weakly connected component $W$ contains 8,123 vertices (users) and 10,911 edges (tweet relations). We cluster $W_{\BC}$ and $W_{\CC}$ to 10 clusters at level 2 and then further cluster them to 2 clusters at level 1.

For this data set, we can perform both unsupervised learning  and semi-supervised learning methods.  The results are given in Table~\ref{proprestab:F-M}. From the table,  in terms of the $F$-measure,  semi-supervised method {\Percus} performs better than the other two methods of {\GCN} and {\LL} while within unsupervised methods, {\LL} is better than {\GCN}. In terms of the  modularity metric, {\GCN} performs in general better than the other two methods especially when the size of the training data is small ($\le 60\%$).

\begin{table}[h]
\begin{center}
\begin{tabular}{c|l|c|c|c|c|c|c|c|c|c|c}
\hline 
\hline
$\mathcal{F}$ & Trains (\%) & 0 (USL) &  10 & 20 & 30 & 40 & 50 & 60 & 70 & 80 & 90  \\
\hline
\multirow{2}{*}{\GCN} &
Level 2 (10)& 0.15& 0.11& 0.15& 0.25& 0.31& 0.39& 0.48& 0.60& 0.73& 0.86\\
& Level 1 (2)& 0.65& 0.55& 0.49& 0.47& 0.47& 0.49& 0.58& 0.64& 0.73& 0.83\\
\hline
\multirow{2}{*}{\LL} &
Level 2 (10)& 0.49& 0.44& 0.43& 0.44& 0.47& 0.52& 0.59& 0.67& 0.76& 0.87\\
& Level 1 (2)& 0.80& 0.77& 0.68& 0.69& 0.65& 0.62& 0.67& 0.67& 0.74& 0.82\\
\hline
\multirow{2}{*}{\Percus} &
Level 2 (10)& N.A.& 0.50& 0.53& 0.59& 0.64& 0.69& 0.75& 0.81& 0.87& 0.94\\
& Level 1 (2)& N.A.& 0.84& 0.86& 0.86& 0.86& 0.86& 0.88& 0.88& 0.90& 0.91\\
\hline
\hline
$\mathcal{M}$ & Trains (\%) & 0(USL) & 10 & 20 & 30 & 40 & 50 & 60 & 70 & 80 & 90  \\
\hline
\hline
\multirow{2}{*}{\GCN} &
Level 2 (10)& 0.48& 0.44& 0.33& 0.28& 0.15& 0.17& 0.12& 0.08& 0.08& 0.08\\
& Level 1 (2)& 0.23& 0.22& 0.20& 0.19& 0.09& 0.13& 0.08& 0.05& 0.06& 0.08\\
\hline
\multirow{2}{*}{\LL} &
Level 2 (10)& 0.09& 0.07& 0.06& 0.05& 0.04& 0.06& 0.04& 0.04& 0.07& 0.09\\
& Level 1 (2)& 0.08& 0.04& 0.03& 0.02& 0.01& 0.04& 0.03& 0.03& 0.05& 0.07\\
\hline
\multirow{2}{*}{\Percus} &
Level 2 (10)& N.A. & 0.28& 0.25& 0.17& 0.15& 0.13& 0.11& 0.09& 0.10& 0.09\\
& Level 1 (2)& N.A. & 0.07& 0.04& 0.05& 0.04& 0.05& 0.04& 0.05& 0.07& 0.06\\
\hline
\hline
\end{tabular}
\end{center}
\caption{Proposition Data: the $F$-measure $\mathcal{F}$ and modularity metric $\mathcal{M}$ using {\GCN, \LL}, and {\Percus} for given 0\% (USL method), 10\%, 20\%, \ldots, 90\% training data at level 2 ($k_2=10$ clusters) and level 1 ($k_1=2$ clusters), respectively. All results are average over 30 trials.}
\label{proprestab:F-M}
\end{table}

\subsection{The Wiki--votes data set}\label{wikisect}

The Wiki--votes data set \cite{leskovec2010signed, leskovec2010predicting} is available from the Stanford large network data set collection at 
 https://snap.stan-ford.edu/data/wiki-Vote.html. Per this website: ``A small part of Wikipedia contributors are administrators, who are users with access to additional technical features that aid in maintenance. In order for a user to become an administrator a Request for adminship (RfA) is issued and the Wikipedia community via a public discussion or a vote decides who to promote to adminship. Using the latest complete dump of Wikipedia page edit history (from January 3 2008) we extracted all administrator elections and vote history data. This gave us 2,794 elections with 103,663 total votes and 7,066 users participating in the elections (either casting a vote or being voted on). Out of these 1,235 elections resulted in a successful promotion, while 1,559 elections did not result in the promotion. About half of the votes in the dataset are by existing admins, while the other half comes from ordinary Wikipedia users. 
The network contains all the Wikipedia voting data from the inception of Wikipedia till January 2008. Nodes in the network represent wikipedia users and a directed edge from node $i$ to node $j$ represents that user $i$ voted on user $j$.''

The graph from the Wiki--votes data set has 7,115 vertices and 103,689 edges. It has 24 weakly connected components and 5,816 strongly connected components. The largest weakly connected component  contains 7,066 vertices and others contains at most 3 vertices and can be viewed as singletons. Hence, we only consider the largest weakly connected component.   

The Wiki--votes data set is unlabelled. Therefore, only the
{\GCN} and {\LL}  algorithms  can be used, and the performance can be only be tested using the modularity metric.  We cluster all vertices to $k_2$ clusters at level 2 and then further cluster them to $k_1$ clusters at level 1.  We choose $k_2$ ranging from $4$ to $11$ and $k_1$ ranging from $2$ to $8$ with step size 1. For each possible pair $(k_1,k_2)$, we compute the modularity metric from averaging over 30 trials. In terms of modularity, we found that the best modularity (see  Table~\ref{wikirestab})  is  $(\mathcal{M}_2,\mathcal{M}_1)=(0.040,0.037)$ with respect to  $(k_2,k_1)=(6,3)$ for {\GCN}, where $\mathcal{M}_2, \mathcal{M}_1$ are the modularity at level 2 and level 1, respectively.  While it is $(\mathcal{M}_2,\mathcal{M}_1)=(0.078,0.074)$ with respect to  $(k_2,k_1)=(4,3)$ for {\LL}.

\begin{table}[h]
\begin{center}
\begin{scriptsize}
\begin{tabular}{c|c|c|c|c|c|c|c|c}
\hline 
\hline
\multicolumn{9}{c}{\GCN $(\mathcal{M}_2,\mathcal{M}_1)$}  \\
\hline
$(k_1,k_2)$ & 4 & 5 & 6  & 7 & 8 & 9 & 10 & 11  \\
\hline
2 & (0.035,0.031)& (0.036,0.029)& (0.032,0.030)& (0.035,0.024)& (0.032,0.031)& (0.035,0.022)& (0.032,0.025)& (0.029,0.016)\\
3 & (0.031,0.032)& (0.036,0.036)& (0.040,0.037)& (0.033,0.031)& (0.033,0.034)& (0.032,0.029)& (0.031,0.035)& (0.027,0.030)\\
4 & & (0.032,0.032)& (0.034,0.034)& (0.033,0.034)& (0.030,0.033)& (0.033,0.033)& (0.031,0.036)& (0.028,0.032)\\
5 & & & (0.032,0.032)& (0.035,0.035)& (0.032,0.036)& (0.031,0.032)& (0.030,0.032)& (0.029,0.033)\\
6 & & & & (0.035,0.035)& (0.033,0.034)& (0.031,0.032)& (0.029,0.031)& (0.029,0.032)\\
7 & & & & & (0.032,0.033)& (0.032,0.033)& (0.029,0.031)& (0.029,0.032)\\
8 & & & & & & (0.030,0.030)& (0.030,0.031)& (0.030,0.032)\\
\hline
\hline
\multicolumn{9}{c}{\LL $(\mathcal{M}_2,\mathcal{M}_1)$}  \\
\hline
$(k_1,k_2)$ & 4 & 5 & 6  & 7 & 8 & 9 & 10 & 11  \\
\hline
2 & (0.078,0.058)& (0.056,0.052)& (0.048,0.048)& (0.046,0.045)& (0.041,0.036)& (0.036,0.038)& (0.032,0.042)& (0.028,0.042)\\
3 & (0.078,0.074)& (0.056,0.060)& (0.047,0.050)& (0.046,0.048)& (0.040,0.044)& (0.036,0.040)& (0.032,0.042)& (0.030,0.040)\\
4 & & (0.056,0.061)& (0.048,0.049)& (0.045,0.046)& (0.040,0.046)& (0.035,0.040)& (0.033,0.045)& (0.031,0.039)\\
5 & & & (0.047,0.049)& (0.046,0.050)& (0.040,0.044)& (0.037,0.043)& (0.032,0.041)& (0.029,0.041)\\
6 & & & & (0.046,0.049)& (0.040,0.041)& (0.036,0.041)& (0.032,0.036)& (0.030,0.034)\\
7 & & & & & (0.040,0.041)& (0.036,0.038)& (0.033,0.036)& (0.029,0.033)\\
8 & & & & & & (0.036,0.039)& (0.032,0.035)& (0.030,0.033)\\
\hline
\hline
\end{tabular}
\end{scriptsize}
\end{center}
\caption{Wiki--vote Data: the modularity  metric pair $(\mathcal{M}_2, \mathcal{M}_1)$  for the Wiki-Vote data set using GCN and LL unsupervised clustering  ($\mathcal{M}_2$ is the modularity at level 2 and $\mathcal{M}_1$ is the modularity at level 1) for each $(k_1,k_2)\in\{2,3,\ldots,8\}\times \{4,5,\ldots,11\}$. All results are average over 30 trials.}
\label{wikirestab}
\end{table}

\bhag{Approximation and analysis of functions}\label{harmanalsect}

To the best of our knowledge, representation of functions 
has been accomplished typically by a spectral decomposition of a graph Laplacian. We cite, for example, \cite{hammond} for undirected graphs, which follows ideas developed in \cite{fasttour}, and \cite{gidelew2014topics} for directed graphs. A more analytical treatment is typically based on embedding the digraph into a manifold \cite{mousazadeh2015embedding, tauberian}. Our paper \cite{treepap} gave a completely different approach that is free of any spectral decomposition, except for spectral clustering methods if used. 
The purpose of this section is to extend the results in \cite{treepap} to the case of digraphs, represented by a pair of trees.

In this section, we assume that the trees corresponding to the digraph have been constructed using  appropriate clustering algorithms. In greater abstraction, we will assume that the trees are both infinite. In practical terms, this allows us to add data to the digraph. We will 
then describe harmonic analysis on this infinite digraph. In Sub--Section~\ref{reviewsect}, we will review the relevant ideas from \cite{treepap}. They will be used in Sub--Section~\ref{diharmsect} to
describe a very general harmonic analysis and function approximation paradigm.

\subsection{Tree polynomials}\label{reviewsect}

Fundamental to our analysis is the notion of a filtration, as defined below.
\begin{definition}\label{filtrationdef}
 A \textbf{weighted tree} is a triplet $(V,E,w)$, where $(V,E)$ is a tree and $w : V\to (0,\infty)$ is a weight function. Let $\XX=(V,E,w)$ be a weighted tree, and $v\in V$.   The tree with root $v$ and leaves given by the children of $v$ will be denoted by $\XX(v)$. The tree $\XX(v)$ is called the \textbf{local (or conditional) filtration} at $v$, if the set $\mathcal{L}_v$ of the children of $v$ contains at least $2$ elements and $\sum_{u\in \mathcal{L}_v} w(u)=w(v)$. The weighted tree $\XX$ is called a \textbf{filtration}, if $w(\rr)=1$; and for each non--leaf vertex $v\in V$, $\XX(v)$ is a local filtration at $v$. 
\end{definition}

It is possible to represent the vertices of a filtration as subintervals of $[0,1)$. Let $\XX=(V,E,w)$ be a filtration. We associate the root $\rr$ with the unit interval $[0,1)$. If the children of $\rr$ are $v_1,\cdots,v_{M}$, we associate with each $v_j$ the interval 
$$
\left[\sum_{k=1}^{j-1}w(v_k), \sum_{k=1}^{j}w(v_k)\right).
$$
Our assumption that $\sum_{k=1}^{M}w(v_k)=w(\rr)=1$ implies that these intervals constitute a partition of $[0,1)$.  In a recursive manner, if the  interval associated with a vertex $v$ is $[a,a+w(v))$, and the children of $v$ are $u_1,\cdots,u_K$, then we associate with each $u_j,\ j=1,\cdots,K,$ the interval 
$$
\left[a+\sum_{k=1}^{j-1} w(u_k), a+\sum_{k=1}^j w(u_k)\right).
$$
Since $\XX(v)$ is a local filtration at $v$, the intervals associated with the $u_j$'s constitute a partition of the interval 
associated with $v$.  In the sequel, we will refer to a vertex on the tree and the associated interval interchangeably. If $S\subseteq [0,1)$, we denote by $\cfn{S}$ the characteristic function of $S$; i.e., $\cfn{S}(x)=1$ if $x\in S$, and $\cfn{S}(x)=0$ if $x\in [0,1)\setminus S$.

For a local filtration $\XX(v)$ at $v$ associated with interval $[a,a+w(v))$ for some $a$, let $v_0,\cdots, v_{m-1}$, $m\ge 2$, be the children of $v$ in order. We wish to obtain the set of tree polynomials  to have the same span as $\{ \cfn{v_j} : \ j=0,\cdots,m-1\}$. In order to facilitate the construction of a consistently labeled system across the entire filtration, it is convenient to 
substitute $\cfn{v_0}$ by $\cfn{v}$ in the list above.  The resulting system is then defined by \eref{localtreepolydef} as follows. Let $p_j=w(v_j)$ for $j=0,\dots, m-1$ and
\be\label{localpdef}
 P_k:=\sum_{j=0}^{k-1}p_j, \quad I_k:=[a+P_k,a+P_{k+1}), \quad J_k:=[a,a+P_k), \qquad k=1,\cdots,m-1,  \quad P_0:=0, \ P_{m+1}=w(v).
\ee
Define 
\be\label{localtreepolydef}
\begin{aligned}
\phi_0(x)&=\phi_0(v,x)=\cfn{v},\\
 \phi_k(x)& =\phi_k(v,x)=p_k\cfn{J_k}(x)-P_k\cfn{I_k}(x), \qquad x\in [a,a+w(v)),\ k=1,\cdots,m-1.
 \end{aligned}
\ee
Some of the important properties of these tree polynomials are listed in the following proposition (\cite[Proposition~3.1, Proposition~3.2]{treepap}).
\begin{prop}\label{locatreeprop}
Let $\XX(v)$ be a local filtration, and $\phi_j$'s be defined as in \eref{localtreepolydef}. 
Let $\ell, k$ be integers $0\le \ell, k\le m-1$. Then\\
{\rm (a)}  $\phi_0 = \cfn{v}$ and for $k=1,\ldots, m-1$,
\be\label{locphivalues}
\phi_k(x)= \left\{\begin{array}{ll}
p_k, &\mbox{ if $x\in I_\ell$, $1\le \ell \le k-1$,}\\
-P_k, &\mbox{ if $x\in I_k$},\\
0, &\mbox{if $x\in I_\ell$, $\ell \ge k+1$;}
\end{array}\right.
\ee
{\rm (b)} 
\be\label{locphiortho}
\int_v \phi_k(x)\phi_\ell(x)dx=\left\{\begin{array}{ll}
w(v), &\mbox{if $k=\ell=0$,}\\
p_kP_kP_{k+1}, &\mbox{if $k=\ell\ge 1$,}\\
0, & \mbox{otherwise};
\end{array}\right.
\ee
{\rm (c)} 
\be\label{localquad}
\sum_{j=0}^{m-1}p_j\phi_k(a+P_j)\phi_\ell(a+P_j)=\int_v\phi_k(x)\phi_\ell(x)dx,
\ee
and for $j=0,\cdots,m-1$,
\be\label{localinvquad}
1+\sum_{k=1}^{m-1}\frac{\phi_k(a+P_j)\phi_k(a+P_\ell)}{p_kP_kP_{k+1}}=\left\{\begin{array}{ll}
1/p_j, &\mbox{if $j=\ell$,}\\
0, &\mbox{otherwise.}
\end{array}\right.
\ee
{\rm (d)} Let $0\le k\le m-1$ be an integer, and for $x\in [a,a+w(v))$,
$$
\tilde{\phi}_{k+1}(x)=\tilde{\phi}_{k+1}(v;x)=\left\{\begin{array}{ll}
0, &\mbox{if $k=m-1$},\\
(1-P_{k+1})\cfn{J_{k+1}}(x)-P_{k+1}\cfn{v\setminus J_{k+1}}(x), &\mbox{if $k\le m-2$.}
\end{array}\right.
$$
Then 
\be\label{temportho}
\int_v\phi_j(x)\tilde{\phi}_{k+1}(x)dx=0, \qquad j=0,\cdots, k,
\ee
and
\be\label{prelocalspan}
\tilde{\Pi}_k:=\tilde{\Pi}_k(v):=\spans\{\phi_0,\cdots,\phi_k,\tilde{\phi}_{k+1}\}=\spans\{\cfn{I_0},\cdots,\cfn{I_k},\cfn{v\setminus J_{k+1}}\}.
\ee
\end{prop}

The orthogonal system on the whole tree is designed, so that in principle, when restricted to each local filtration, it should reduce to the local system for that filtration. To describe this in detail, we need to introduce the notion of Leave Left Out (LLO) enumeration. At each level $L\ge 1$, let (in this section only) 
$$
C_L=\{v : v \mbox{ is a vertex at level $L$, $v$ is not a left--most child of its parent}\},
$$
$C_0=\{\rr\}$. We associate $M_0=0$ with the root $\rr$ at level $0$. At each level $L\ge 1$, we enumerate the vertices in $C_L$ by $M_{L-1},\cdots,M_L-1$, left to right. This enumeration will be called the Leave Left Out (LLO) enumeration. 

For example, the LLO enumeration of the vertices of the  left tree in Figure~\ref{toydaspict}  is $\{A, B, E, F, H, I\}$. We note that vertices which are 
left--most children of their parents are not numbered in this scheme.

Each vertex in $V\setminus\{\rr\}$ can also be enumerated among the children of its parent. For any $v\in V$, the children of 
$v$ are enumerated left to right starting from $0$, with $0$ associated with the left--most child of $v$.

\begin{definition}\label{globalorthosystemdef}
For $x\in [0,1)$, let $\psi_0(x)=1$; and for an integer $n\ge 1$, let $u$ be the vertex corresponding to $n$ under the LLO enumeration, $v$ be the parent of $u$, and $\ell$ be the enumeration  of $u$ as a child of $v$, so that $\ell \ge 1$. Set
\be\label{psindef} 
\psi_n(x)=\phi_\ell(v,x).
\ee
The symbol $\Pi_n$ denotes the span of $\{\psi_0,\cdots,\psi_n\}$.
\end{definition}

The following proposition \cite[Theorem~4.2]{treepap} summarizes some of the properties of the system $\{\psi_n\}$.

\begin{prop}\label{globaltreeprop}
Let $\XX$ be a filtration, $n, m\ge 0$ be integers, $u=[a',b')$, $u'$ be vertices corresponding to $n$, $m$ respectively in the LLO enumeration, $v=[a,b)$ be the parent of $u$. Then
\be\label{fullorthoeqn}
\int_0^1 \psi_n(x)\psi_m(x)dx =\left\{\begin{array}{ll}
0, &\mbox{ if $n\not=m$,}\\
1, &\mbox{if $n=m=0$,}\\
\aleph_n^{-1}:=\disp\frac{(b'-a')(a'-a)(b'-a)}{(b-a)^2}, &\mbox{ if $n=m\not=0$.}
\end{array}\right.
\ee 
\end{prop}

If $f:[0,1)\to \RR$ is any bounded and integrable function, then we define for integers $k\ge 0$, $n\ge 1$,
\be\label{treefourdef}
\hat{f}(k)=\aleph_k\int_0^1 f(t)\psi_k(t)dt, \qquad s_n(f,x)=\sum_{k=0}^{n-1} \hat{f}(k)\psi_k(x).
\ee
A novelty of our system is that the partial sum operators $\{s_n\}$ themselves are uniformly bounded, in contrast with the classical theory of Fourier series.
Analogous to the summability methods in the theory of Fourier series,
we define a more general version of these operators. 

If $\mathbf{h}=\{h_k\}_{k=0}^\infty$ is any sequence, we define 
\be\label{gensummopdef}
\sigma_n(\mathbf{h},f)=\sum_{k=0}^n h_k \hat{f}(k)\psi_k, \qquad n=0,1,\cdots.
\ee
We emphasize again that since the tree polynomials are piecewise constants, both the quantities $s_n(f)$ and $\sigma_n$ can be computed exactly as discrete sums, even though we find it convenient for theory and exposition purposes to write them as integral operators. Further details on this matter are given in \cite{treepap}.

We end this sub--section by enumerating some relevant properties of the operators $\sigma_n$. In the remainder of this sub--section, let $\|\cdot\|$ denote the uniform norm on $[0,1)$, $\mathcal{F}$ be the closure of $\disp\bigcup_{n=0}^\infty \Pi_n$ in this norm, and $E_n(f):=\inf_{P\in\Pi_n}\|f-P\|$ be the degree of approximation of a function $f$ defined on $[0,1)$.

\begin{theorem}\label{genfourtheo}
Let $\mathbf{h}=\{h_k\}_{k=0}^\infty$ be a sequence of real numbers with
\be\label{hseqcond}
\mathcal{V}(\mathbf{h})=\sup_{n\ge 0} \left\{|h_n|
+\sum_{k=0}^{n-1} |h_{k+1}-h_k|\right\} <\infty, \qquad \lim_{N\to\infty} h_N=0.
\ee
We have 
\be\label{genseqfourbd}
\|\sum_{k=0}^n h_k\hat{f}(k)\psi_k\| \le 3\mathcal{V}(\mathbf{h})\|f\|, \qquad f\in \mathcal{F},\quad n=0,1,\cdots.
\ee
In addition, if $h_k=1$ for $0\le k\le n$ for some integer $n\ge 0$, then for $f\in \mathcal{F}$, 
\be\label{gengoodapprox}
E_n(f)\le \|f-\sigma_n(\mathbf{h},f)\|\le (1+3\mathcal{V}(\mathbf{h}))E_n(f).
\ee
\end{theorem}

\subsection{Harmonic analysis on digraphs}\label{diharmsect}
 
In the case of a filtration, there is a one--to--one correspondence between functions on $[0,1)$ and functions on the vertices of the filtration, and it is a matter of convenience whether one thinks of the Lebesgue measure or a discrete measure defined on these vertices. 
When we represent a digraph using a pair of filtrations, this is no longer the case.  It is seen already in the toy example in the introduction (cf. Figure~\ref{toydaspict}), that the points on the digraph correspond only to 25 out of the 625 sub--rectangles of $I^2$. 
In the limiting case, as the number of nodes on the digraph tends to infinity, it may or may not happen that the set of points on the square that correspond to the nodes on the digraph is a set of two dimensional Lebesgue measure $0$. The structure of this set is hard to stipulate mathematically, since it depends upon the exact construction of filtrations, which in turn depends upon the particular digraph in question.
Therefore, harmonic analysis and function approximation on digraphs in our paradigm is more delicate than the multivariate analogues of the univariate analysis in classical situations such as trigonometric series, splines, etc. It is not possible to give universal constructions in our paradigm. In this section, we outline an abstract, data--driven  theory. 
 In this section, we will use standard multivariate notation.
 
We  denote the two filtrations corresponding to the given digraph by $\XX_1$, $\XX_2$ respectively, and the set of points on the square $I^2$ that correspond to the nodes on the digraph by $\GG$. Thus, $\x=(x_1,x_2)\in\GG$ if and only if $x_1$ corresponds to the same vertex on the digraph as a node on the filtration $\XX_1$ as the vertex corresponding to the point $x_2$ as a node on the filtration $\XX_2$ (see Figure~\ref{toygraphmatpict} for an example.) If $w_j$ is the weight associated with $x_j$ as a node on the filtration $\XX_j$, then we associate the weight $w_1w_2$ with the point $\x\in\GG$. The resulting measure will be denoted by $\nu^*$. We note that the measure 
$\nu^*$ is a probability measure, but may well be singular with respect to the two dimensional Lebesgue measure on $I^2$. We will denote the uniform norm on $\GG$ by $\|\cdot\|$.

With an abuse of notation, we denote the \textbf{orthonormalized} system of tree polynomials on these by $\psi_{n,1}$, $\psi_{n,2}$ respectively. Naturally, the tensor product tree polynomials
\be\label{tenstreepoly}
\psi_{\k}(\x)=\psi_{k_1,1}(x_1)\psi_{k_2,2}(x_2), \qquad \k=(k_1,k_2)\in \ZZ_+^2, \ \x=(x_1,x_2)\in I^2,
\ee
are an orthonormal basis for square integrable functions on $I^2$. However, many of these are possibly equal to $0$ when restricted to $\GG$. Let $\Omega=\{\k : \psi_\k|_\GG\not\equiv 0\}$. Since the tree polynomials $\psi_\k$ are constants on rectangles in $I^2$, it is clear that
\be\label{tensortho}
\int_\GG \psi_\k\psi_\m d\nu^* =\int_{I^2}\psi_\k(\x)\psi_\m(\x) d\x, \qquad \k,\m\in \Omega.
\ee 
Thus, $\{\psi_\k\}_{\k\in\Omega}$ is an orthonormal basis for $L^2(\GG,\nu^*)$. The closure of the set $\mathsf{span}\ \{\psi_\k : \k\in\Omega\}$ in the uniform norm of $\GG$ will be denoted by $\mathcal{F}$, abusing again the notation from the case of single filtration.

\begin{figure}[h]
\begin{center}
\includegraphics[width=0.3\textwidth,height=0.27\textwidth]{toygraphasmat_rgb.pdf}
\includegraphics[width=0.32\textwidth,height=0.27\textwidth]{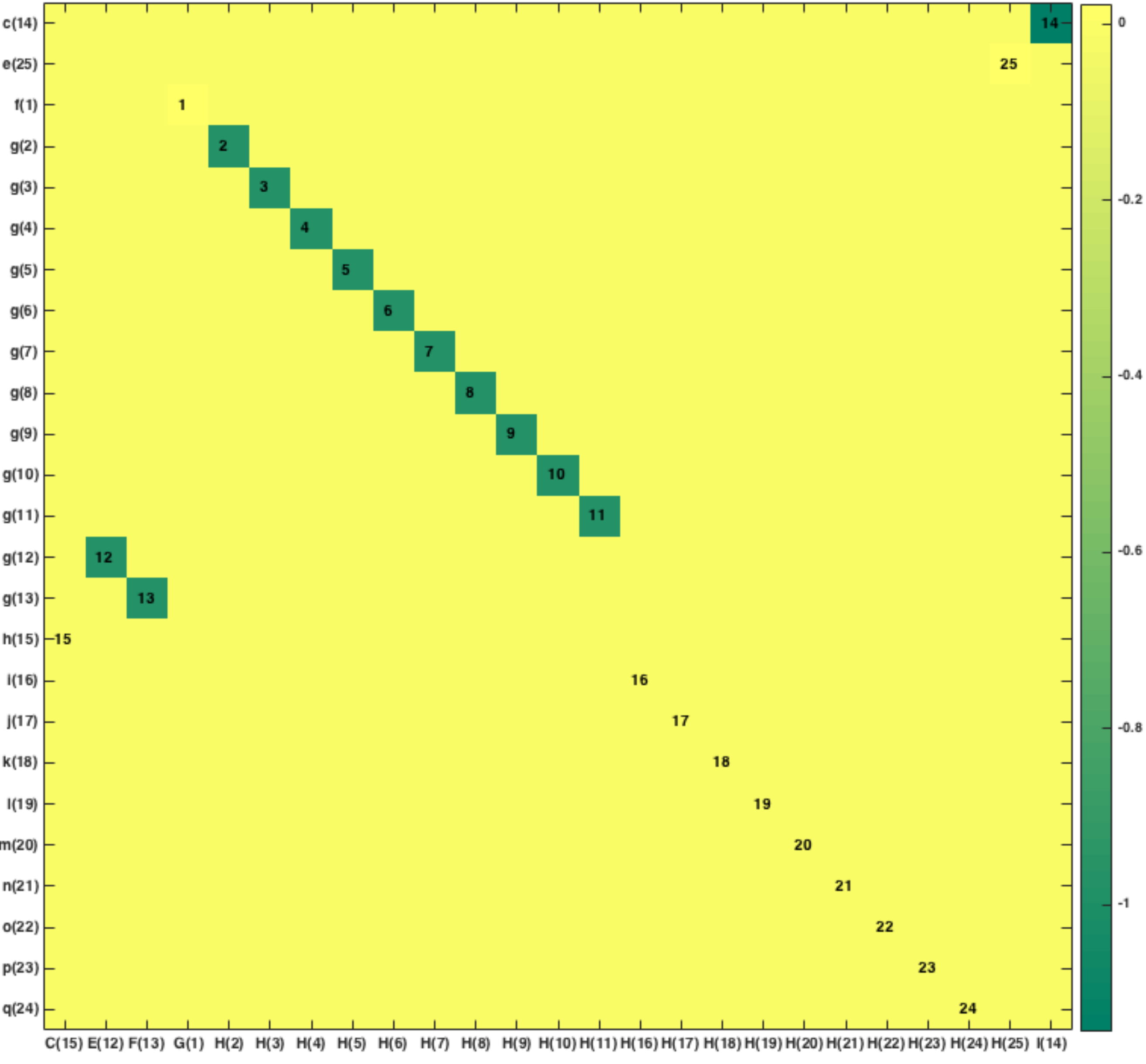}
\includegraphics[width=0.32\textwidth,height=0.27\textwidth]{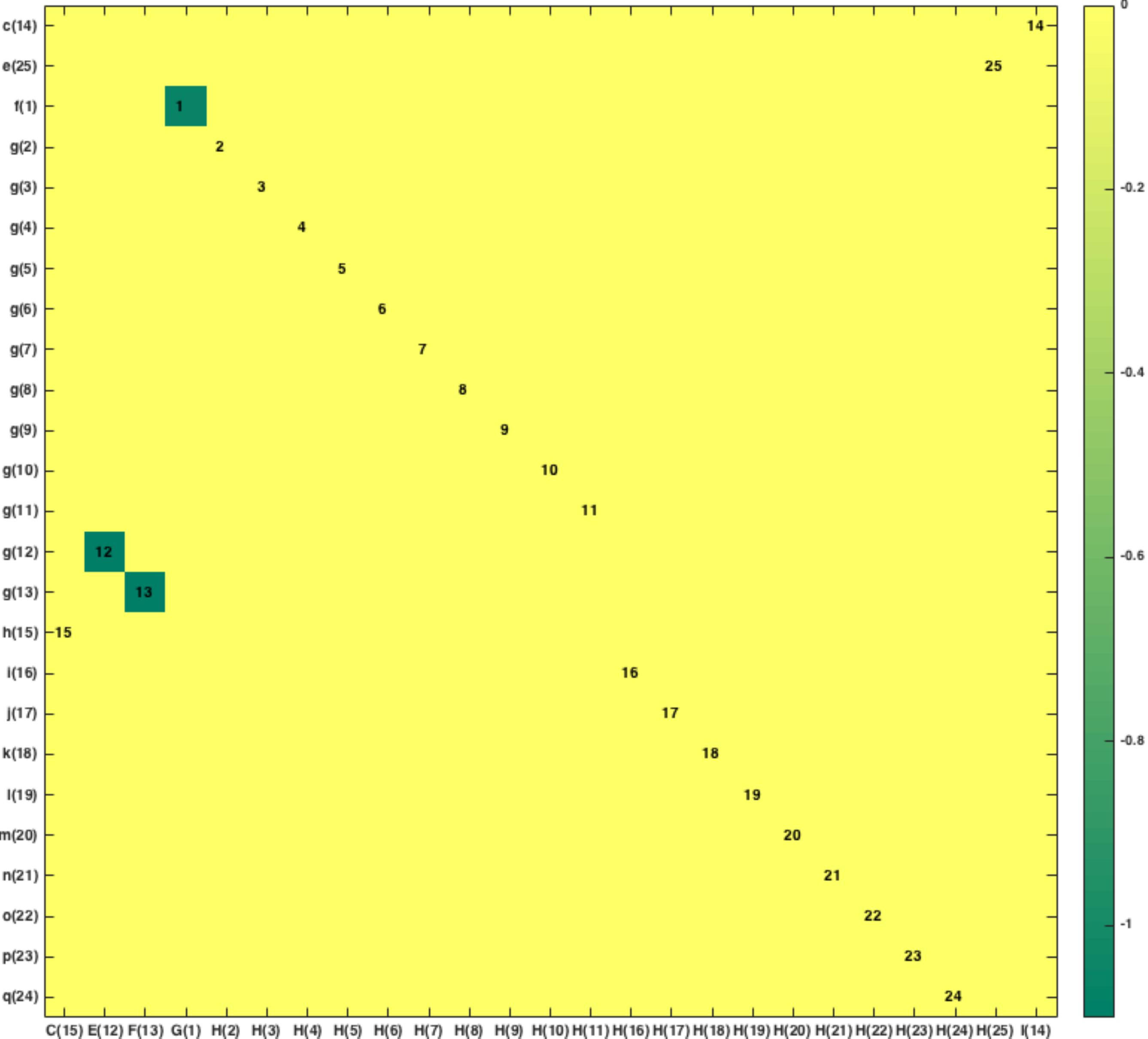}
\end{center}
\caption{The localization of the orthogonal system on $W$ is demonstrated by considering the tensor product tree polynomials for the toy example (see Figure~\ref{toydaspict}). There are $6$ polynomials $\psi_{k_1,1}, k_1=1,\ldots,6$ from $\mathcal{T}_{W_{\BC}}$ and $14$ polynomials $\psi_{k_2,2}, k_2=1,\ldots,14$ from $\mathcal{T}_{W_{\CC}}$. Left: tree vertex position of the tree on a square $[0,1]\times[0,1]$. Middle: tensor product polynomial on the tree w.r.t. $\psi_{2,1}(x)\psi_{6,2}(y)$. Right: tensor product polynomial on the tree w.r.t. $\psi_{3,1}(x)\psi_{6,2}(y)$.}
\label{treepolypict}
\end{figure}

We note that each of the tree polynomials is a piecewise constant function. The localization of these polynomials is illustrated in Figure~\ref{treepolypict} in the context of the toy example in the introduction, referring to the trees in Figure~\ref{toydaspict}.

Since each $\psi_\k$ is bounded, we may define the Fourier coefficients, respectively, partial sums of $f\in L^1(\GG,\nu^*)\cap L^\infty(\GG,\nu^*)$ by
\be\label{tesnfourdef}
\hat{f}(\k)=\int_\GG f\psi_\k d\nu^*, \quad s_\m(f,\x)=\sum_{k_1=0}^{m_1}\sum_{k_2=0}^{m_2} \hat{f}(\k)\psi_\k(\x), \qquad \k, \m\in\ZZ^2.
\ee
It is understood here that $\hat{f}(\k)=0$ if $\k\in\ZZ_+^2\setminus\Omega$, so that the only nonzero summands in the definition of the partial sum correspond to $\k\in\Omega$. 

For a (bi)--sequence $h$, we define
\be\label{tensdiffdef}
\Delta_1h(\k)=h(k_1+1,k_2)-h(k_1,k_2),\quad \Delta_2h(\k)=h(k_1,k_2+1)-h(k_1,k_2), \quad \bs\Delta h=\Delta_1\Delta_2h=
\Delta_2\Delta_1h,
\ee
and
\be\label{hardykrausdef}
\mathcal{V}(h)=\sup_{\k\in\ZZ_+^2} |h(\k)| +\sup_{k_1\in\ZZ_+}\sum_{k_2=0}^\infty|\Delta_2h(\k)| +\sup_{k_2\in\ZZ_+}\sum_{k_1=0}^\infty|\Delta_1h(\k)| +\sum_{\k\in\ZZ_+^2}|\bs\Delta h(\k)|.
\ee
Denoting by $E_1h(\k)=h(k_1+1,k_2)$, $E_2h(\k)=h(k_1,k_2+1)$, we note the following identity for future use. If $h_1$, $h_2$ are sequences as above, then
$\Delta_j(h_1h_2)=\Delta_j(h_1)E_j(h_2)
+h_1\Delta_j(h_2)$, $j=1,2$, and
\be\label{leibnitz}
 \bs\Delta(h_1h_2)=\bs\Delta(h_1)E_1(E_2h_2)
+\Delta_1h_1E_1(\Delta_2h_2)+ \Delta_1h_2E_2(\Delta_1h_2)+
h_1\bs\Delta(h_2).
\ee

In the sequel, $a\lesssim b$ denotes that $a\le c b$ for a generic constant $c$ that does not depend upon the target function and other obvious variables. The value of these generic constants may change at different occurrences, even within a single formula. By $a\sim b$ we mean $a\lesssim b$ and $b\lesssim a$. In particular,
\be\label{leibineq}
\mathcal{V}(h_1h_2)\lesssim \mathcal{V}(h_1)\mathcal{V}(h_2).
\ee

Using Theorem~\ref{genfourtheo}, it is not difficult to prove the following.
\begin{theorem}\label{tensfourtheo}
Let $h=\{h(\k)\}_{\k\in\ZZ_+^2}$ be a finitely supported (bi)--sequence of real numbers. Then for $f\in\mathcal{F}$,
\be\label{tensgenseqfourbd}
\left\|\sum_{\k\in\ZZ_+^2} h(\k)\hat{f}(\k)\psi_\k\right\| \lesssim\mathcal{V}(h)\|f\|.
\ee
\end{theorem}

In the classical theory of multivariate Fourier series, it is customary to define various notions of the degree of the polynomial: spherical, total, coordinatewise, hyperbolic cross, etc. One could do this in the context of tree polynomials on $\GG$ as well, but since the ``frequencies'' are limited to $\Omega$, it is convenient to define a more parsimonious notion  by defining the analysis spaces first and defining the approximation spaces in terms of these. 
\begin{definition}\label{analfilterdef}
A sequence of sequences ${\bf g}=\{g_j :\Omega\to [0,1]\}_{j=0}^\infty$ is called an \textbf{admissible partition of unity} on $\Omega$ if $g_0(\bs 0)=1$, each $g_j$ is supported on a finite set, 
\be\label{partunity}
\sum_{j=0}^\infty g_j(\k)=1, \qquad \k\in\Omega,
\ee
and the following condition is satisfied:
 There exists an integer $m^*=m^*({\bf g})\ge 0$ such that for $j, j'\ge 0$, $|j-j'|>m^*$, $g_j$ and $g_{j'}$ have disjoint supports; i.e., $g_j(\k)g_{j'}(\k)=0$ for all $\k\in \Omega$.
\end{definition}
In the remainder of this section, we will fix an admissible partition ${\bf g}$ of unity. We set $H_n(\k)=\sum_{j=0}^n g_j(\k)$, $\k\in\Omega$, and define the class of multivariate tree polynomials of (${\bf g}$)--degree $\le n$ by
\be\label{tenstrigpolyclass}
\mathbb{P}_n =\spans\{\psi_\k : H_n(\k)>0, \quad \k\in\Omega\}, \qquad n\in\ZZ_+.
\ee
Since each $g_j$ is finitely supported, so is each $H_n$. 

As before, we define the degree of approximation of $f\in L^\infty(\GG,\nu^*)$ by
\be\label{tensdegapproxdef}
E_n(f)=\inf\{\|f-P\| : P\in \mathbb{P}_n\}, \qquad n=0,1,\cdots.
\ee
It is convenient to extend this notation to $n\in\RR$ by setting $E_n(f)=\|f\|$ if $n<0$ and $E_n(f)=E_{\lfloor n\rfloor}(f)$ if $n$ is not an integer.

Next, we define the reconstruction and analysis operators (with an abuse of notation) by
\bea\label{tenssigmataudef}
\sigma_n(f)&=&\sum_{\k\in\Omega}H_n(\k)\hat{f}(\k)\psi_\k,\qquad n=0,1,\cdots, \nonumber\\
\tau_j(f)&=&\left\{\begin{array}{ll}
\sigma_0(f)= \hat{f}(\bs 0),&\mbox{ if $j=0$,}\\
\sigma_j(f)-\sigma_{j-1}(f)=\sum_{\k\in\Omega}g_j(\k)\hat{f}(\k)\psi_\k, &\mbox{if $j=1,2,\cdots$.}
\end{array}\right.
\eea
The following theorem lists some important properties of these operators.
\begin{theorem}\label{tensgoodapproxtheo}
Let $\mathcal{V}(H_n)\lesssim 1$ for all $n\ge 1$, and $f\in \mathcal{F}$. \\
{\rm (a)} We have
\be\label{tensgoodapprox}
E_n(f)\le \|f-\sigma_n(f)\|\lesssim E_{n-m^*}(f), \qquad n=0,1,\cdots.
\ee
{\rm (b)} We have
\be\label{tenslpseries}
f=\sum_{j=0}^\infty \tau_j(f),
\ee
where the sum converges uniformly.\\
{\rm (c)} We have
\be\label{frameprop}
\int_\GG |f(\x)|^2d\nu^*(\x)\sim \sum_{j=0}^\infty \int_\GG |\tau_j(f)(\x)|^2d\nu^*(\x)
\ee
\end{theorem}

We digress to make some comments on the approximation power of our scheme in the context of classification problems. 
The value of the target function $f$ for any leaf is its class label. For simplicity, let us consider a local filtration $\XX(v)$, in which the majority of the leaves have a label $1$, the others have a label $0$.
 Assuming that both classes appear with equal probability, the value of the tree polynomial approximation to $f$ at $v$ is the expected value of the labels of the leaves. 
 To view this as a class label, we need to round it to the nearest integer. This amounts to declaring that the label of the class at the level $v$ is the same as that of the majority of the children of $v$. More generally,  each  cluster $C_i$ in $\{C_1,\ldots,C_M\}$  obtained in the construction of
  the digraph  is assigned the class label $i_{j_0}=\mathrm{argmax}_j |C_i\cap L_j|$ by comparing to the ground-truth classes $\{L_1,\ldots,L_n\}$. Then one defines a confusion matrix $M$ of size $n\times n$ by 
\begin{equation}\label{def:cm}
M_{j,k}:=\sum_{i_{j_0}=j}\frac{|C_{i_{j_0}}\cap L_k|}{|L_k|},\quad j,k=1,\ldots, n,
\end{equation}
where $M_{j,k}$ is the $(j,k)$-entry of the matrix $M$. 
Note that the more $M$ closes to the identity matrix, the better the classification result.

\begin{figure}[htb]
\begin{minipage}{\textwidth}
  \centering
  \begin{minipage}{0.16\textwidth}
  \includegraphics[width=0.99\textwidth,height=0.99\textwidth]{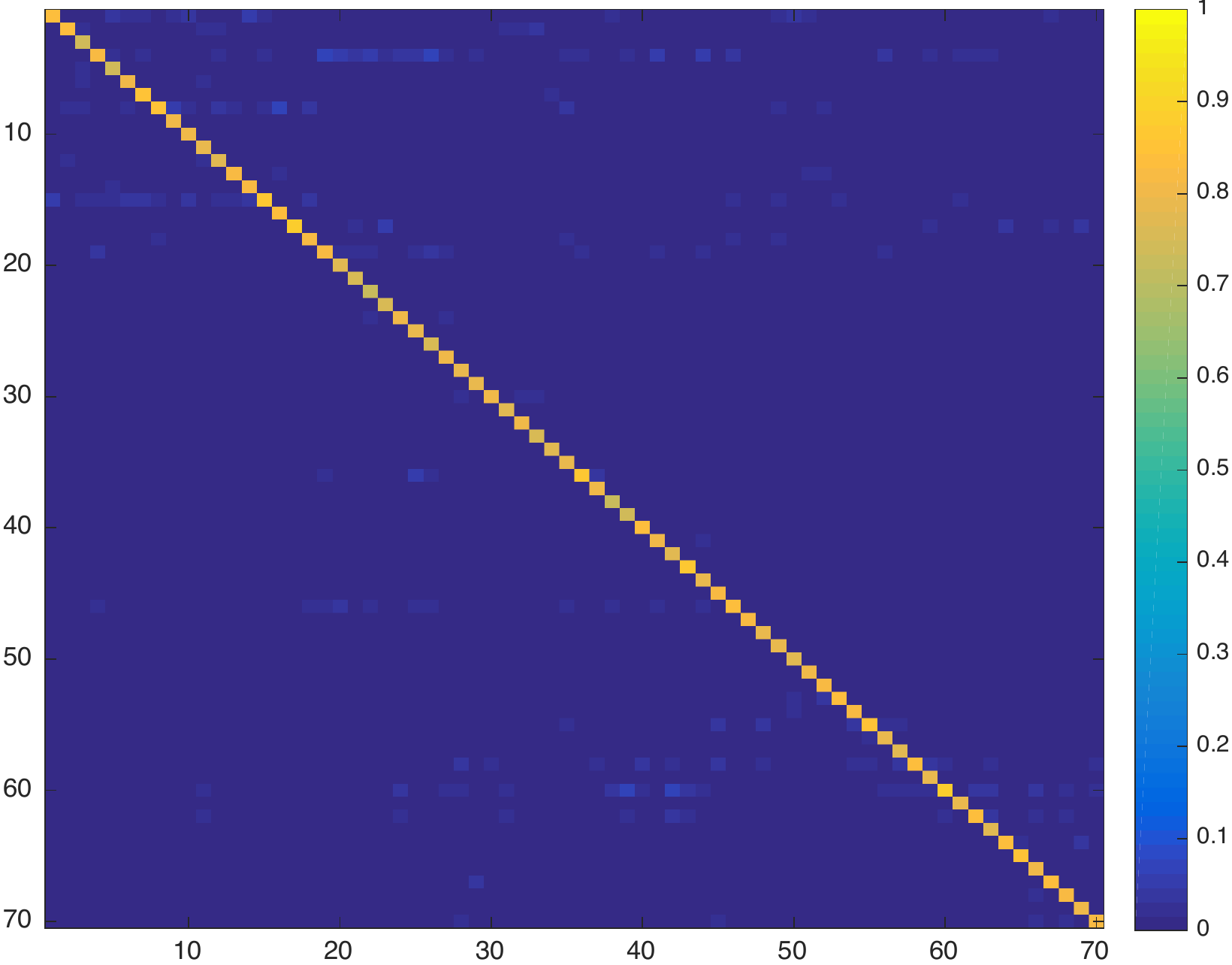}\\[1mm]
  \subcaption{{\GCN}} 
  \end{minipage}
  \begin{minipage}{0.16\textwidth}
  \includegraphics[width=0.99\textwidth,height=0.99\textwidth]{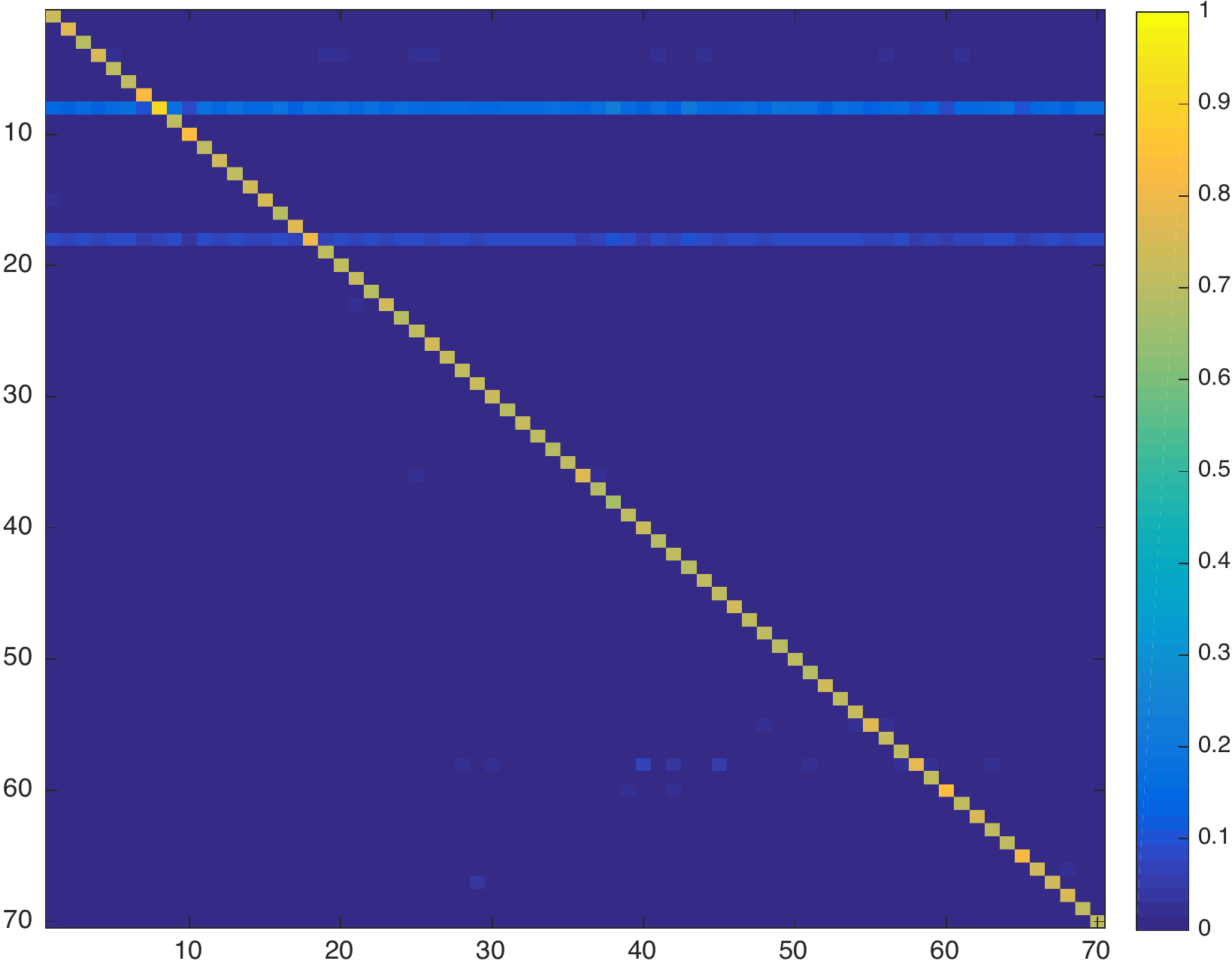}\\[1mm]
  \subcaption{{\LL}}  
  \end{minipage}
  \begin{minipage}{0.16\textwidth}
  \includegraphics[width=0.99\textwidth,height=0.99\textwidth]{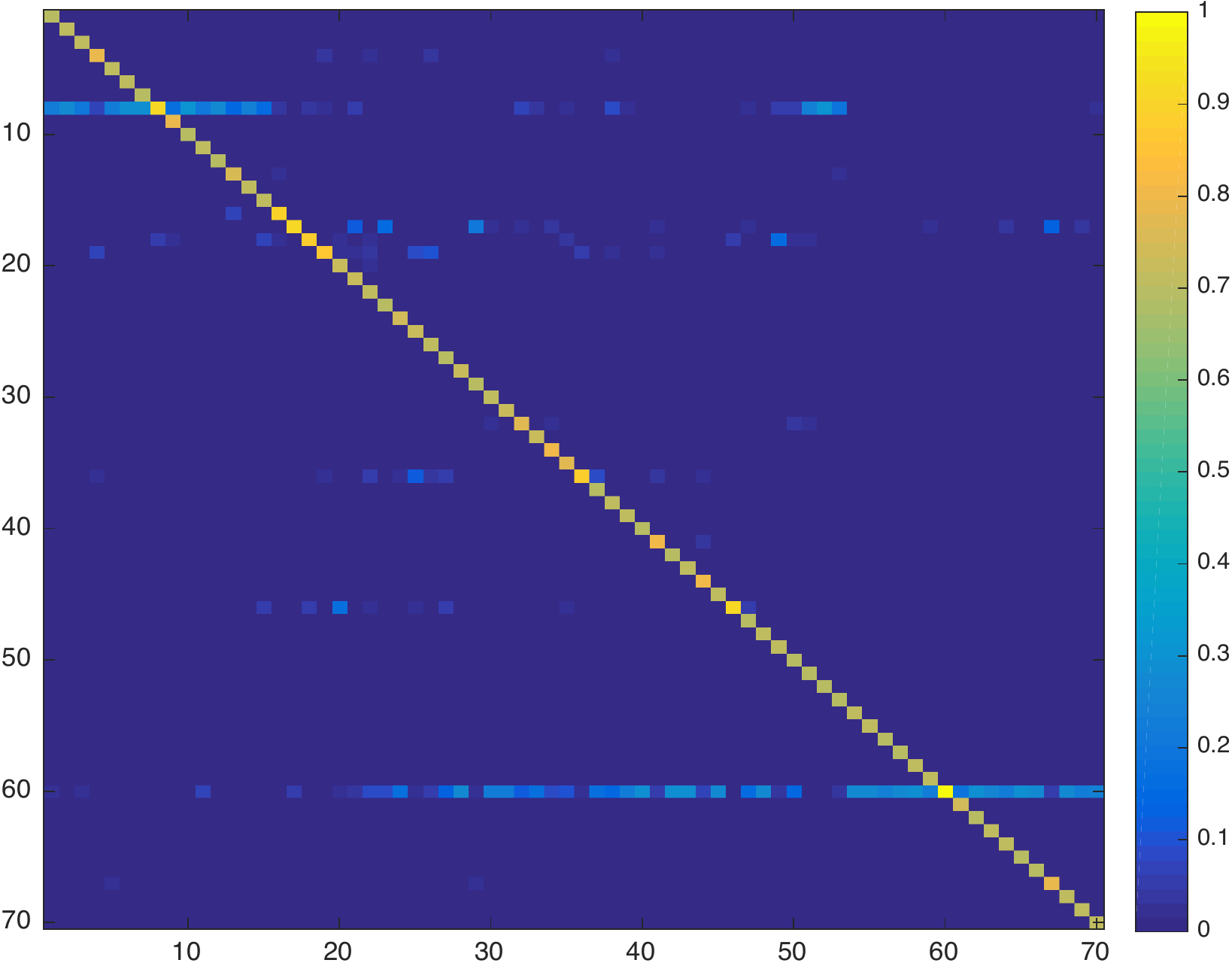}\\[1mm]
  \subcaption{{\Percus}}  
  \end{minipage}
  \begin{minipage}{0.16\textwidth}
  \includegraphics[width=0.99\textwidth,height=0.99\textwidth]{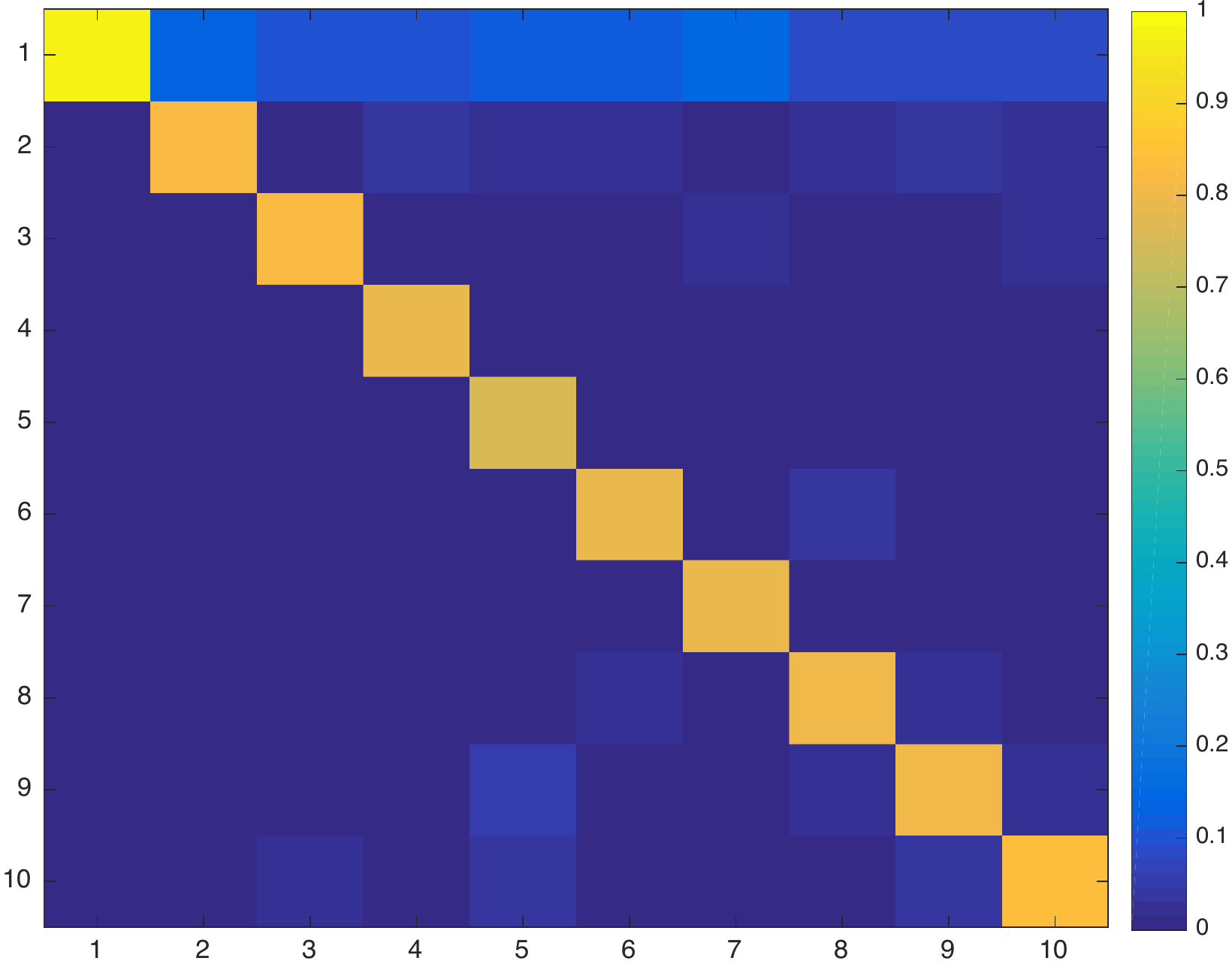}\\[1mm]
  \subcaption{{\GCN}}
  \end{minipage}
  \begin{minipage}{0.16\textwidth}
  \includegraphics[width=0.99\textwidth,height=0.99\textwidth]{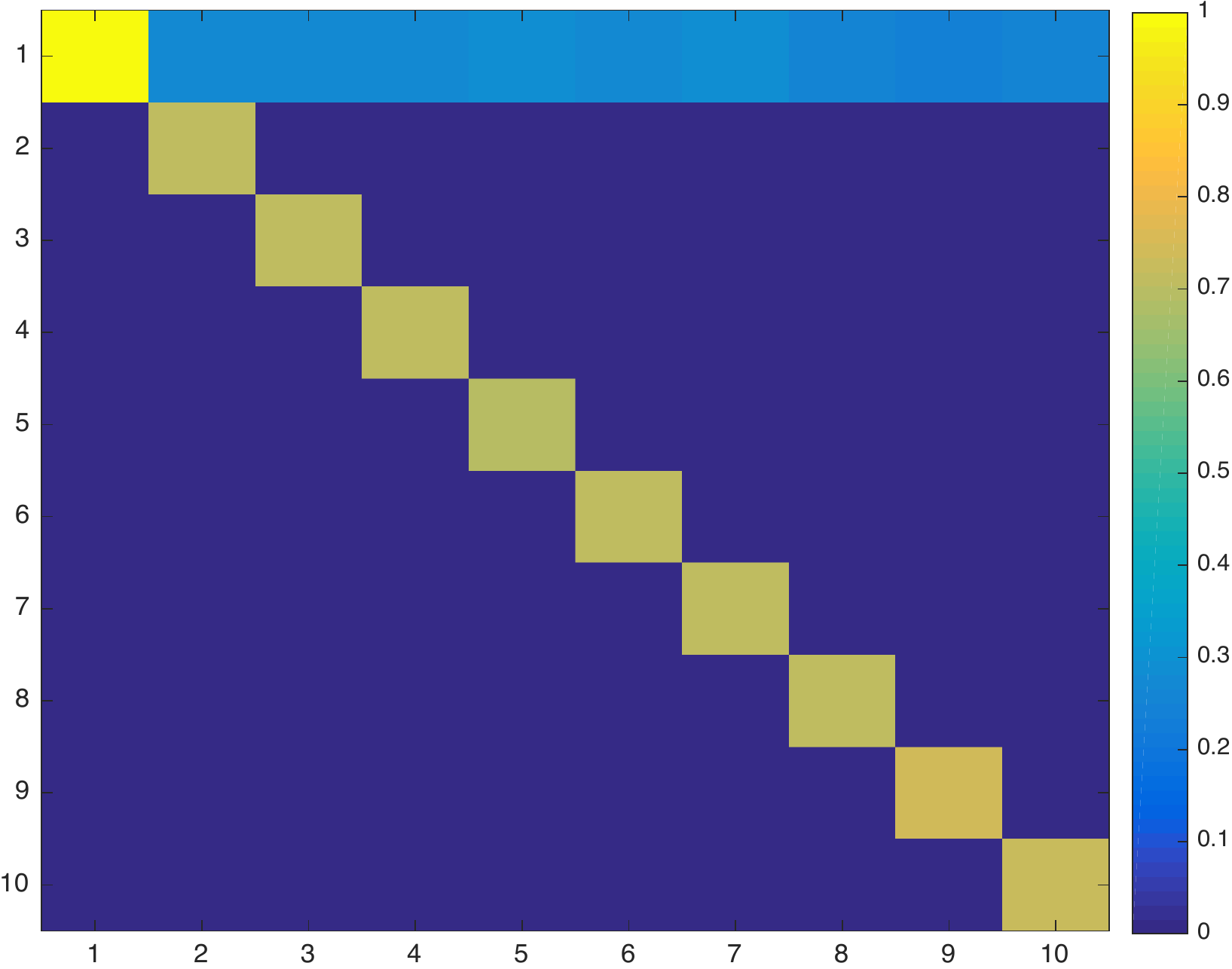}\\[1mm]
  \subcaption{{\LL}}  
  \end{minipage}
  \begin{minipage}{0.16\textwidth}
  \includegraphics[width=0.99\textwidth,height=0.99\textwidth]{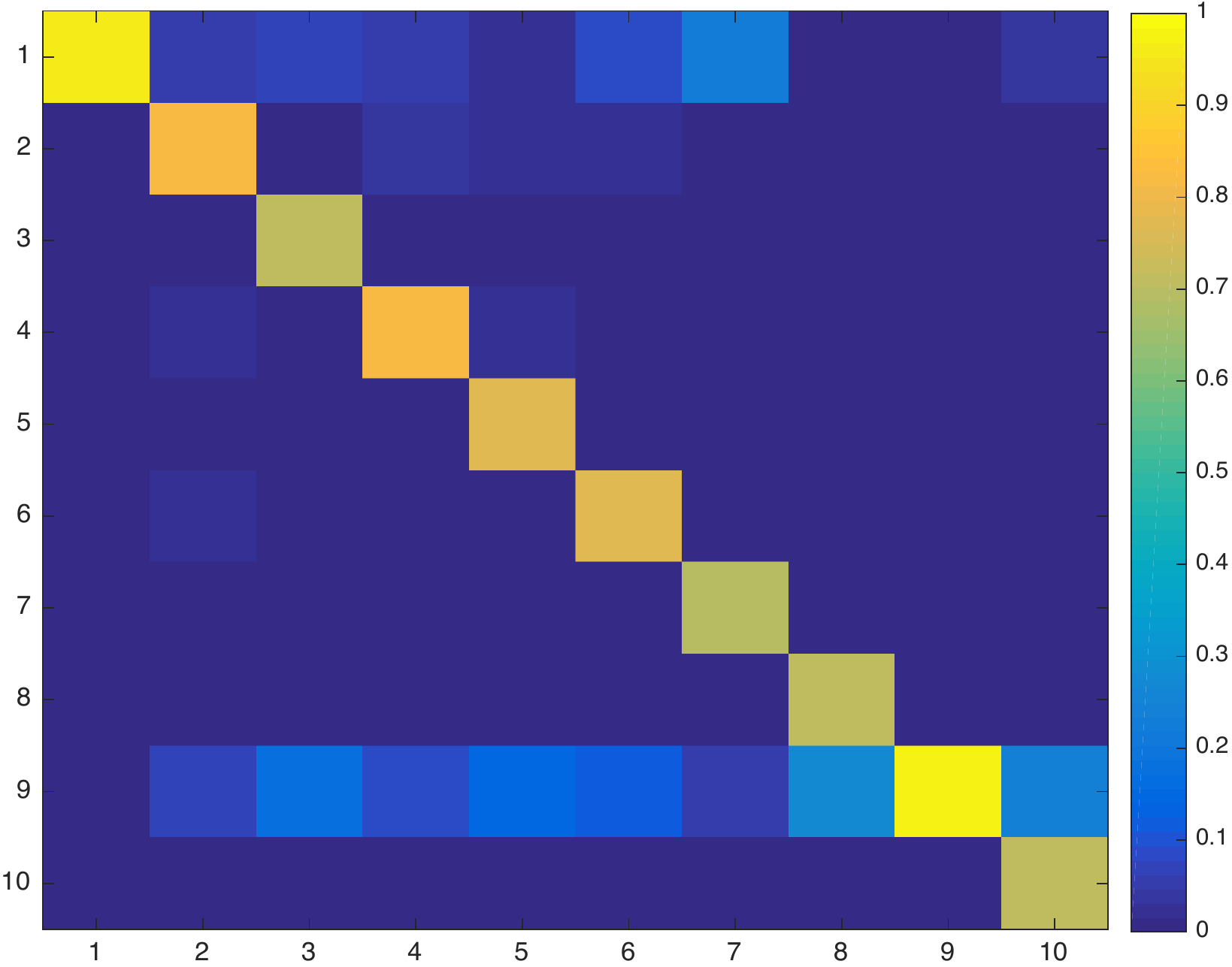}\\[1mm]
  \subcaption{{\Percus}}  
  \end{minipage}
\end{minipage}
\caption{Confusion matrices (images) for CORA dataset. Left 3 (Level 2: 70 classes). Right 3 (Level 1: 10 clases)}
\label{fig:CM_CORA}
\end{figure}

\begin{figure}[htb]
\begin{minipage}{\textwidth}
  \centering
  \begin{minipage}{0.16\textwidth}
  \includegraphics[width=0.99\textwidth,height=0.99\textwidth]{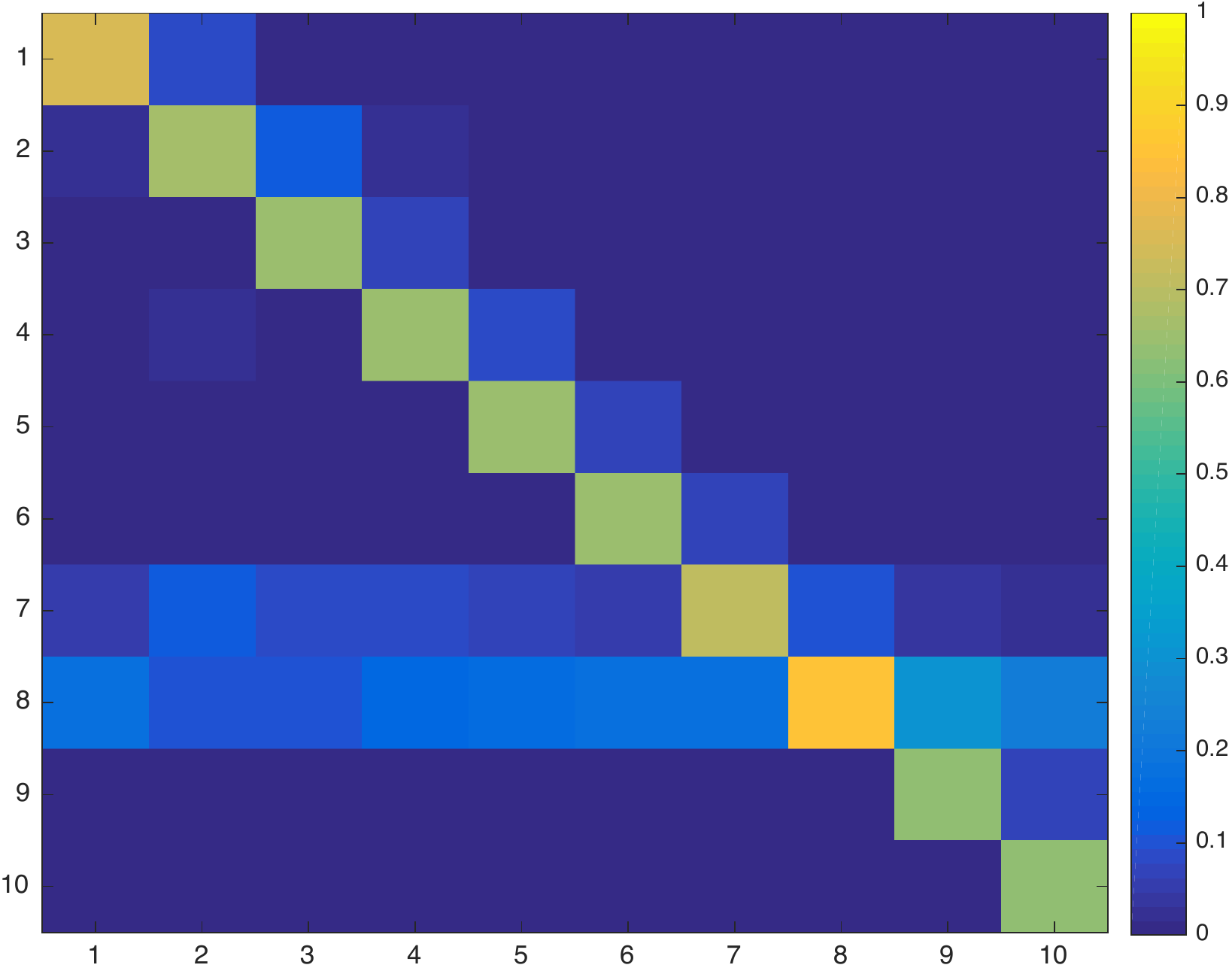}\\[1mm]
  \subcaption{{\GCN}}  
  \end{minipage}
  \begin{minipage}{0.16\textwidth}
  \includegraphics[width=0.99\textwidth,height=0.99\textwidth]{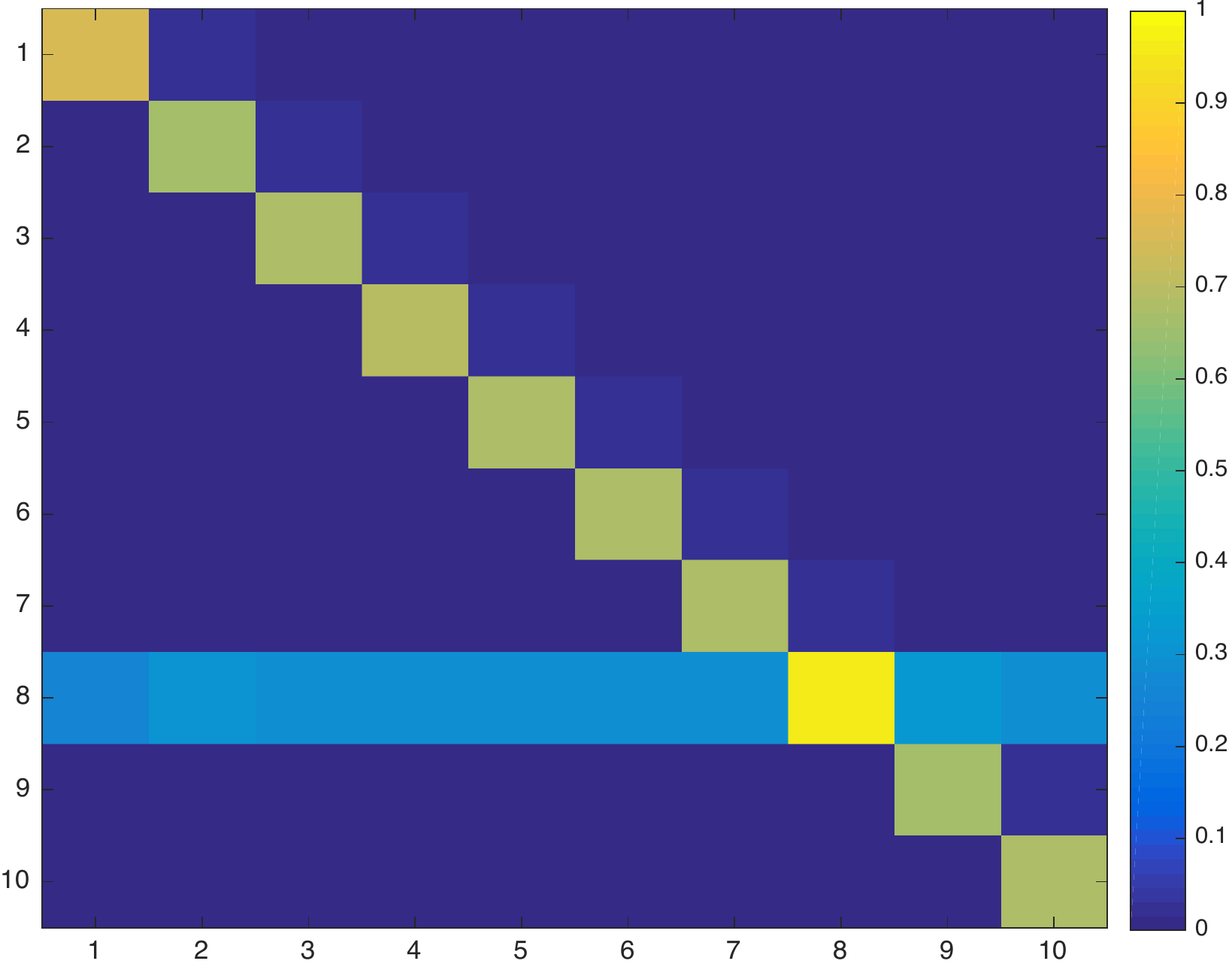}\\[1mm]
  \subcaption{{\LL}} 
  \end{minipage}
  \begin{minipage}{0.16\textwidth}
  \includegraphics[width=0.99\textwidth,height=0.99\textwidth]{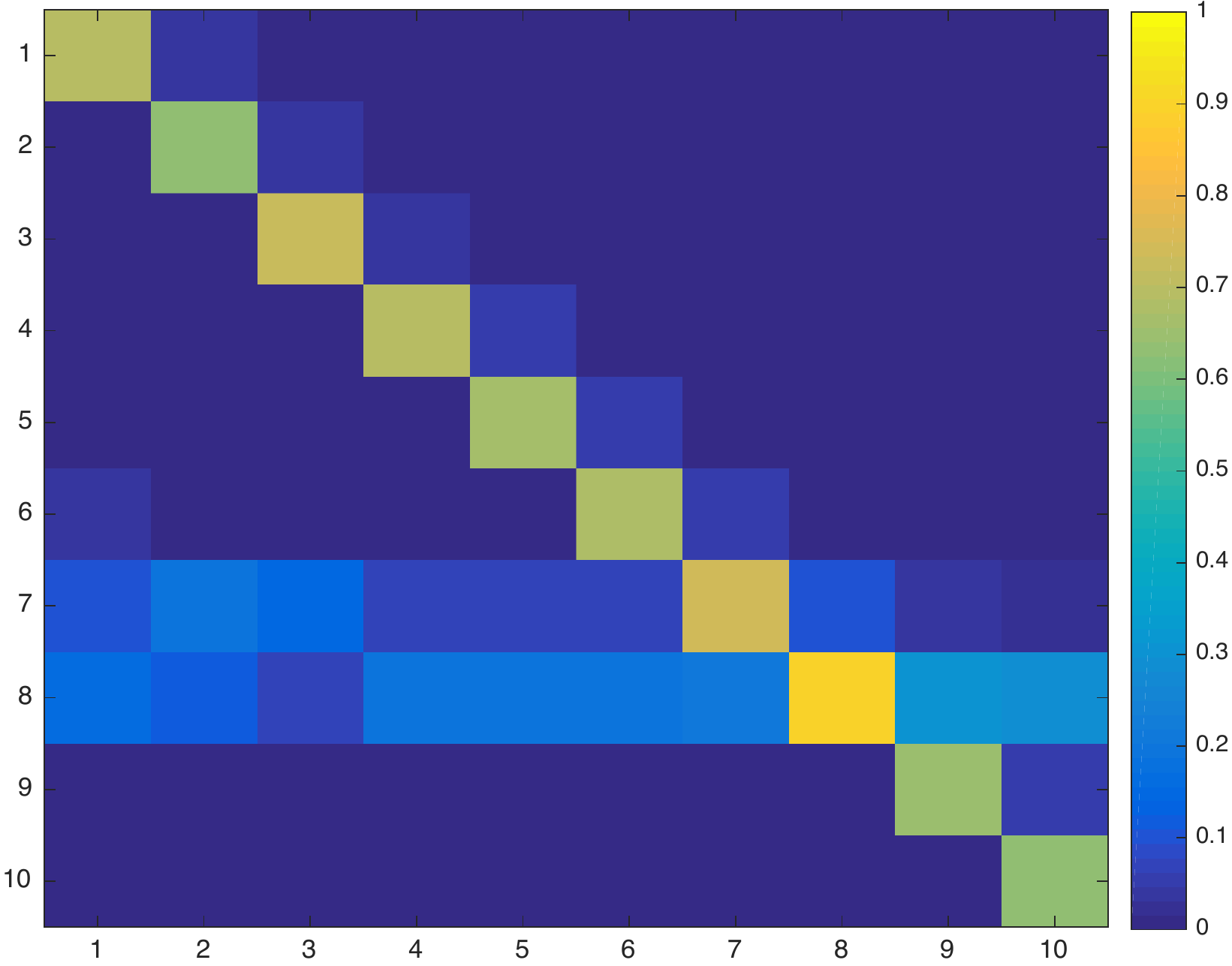}\\[1mm]
  \subcaption{{\Percus}} 
  \end{minipage}
  \begin{minipage}{0.16\textwidth}
  \includegraphics[width=0.99\textwidth,height=0.99\textwidth]{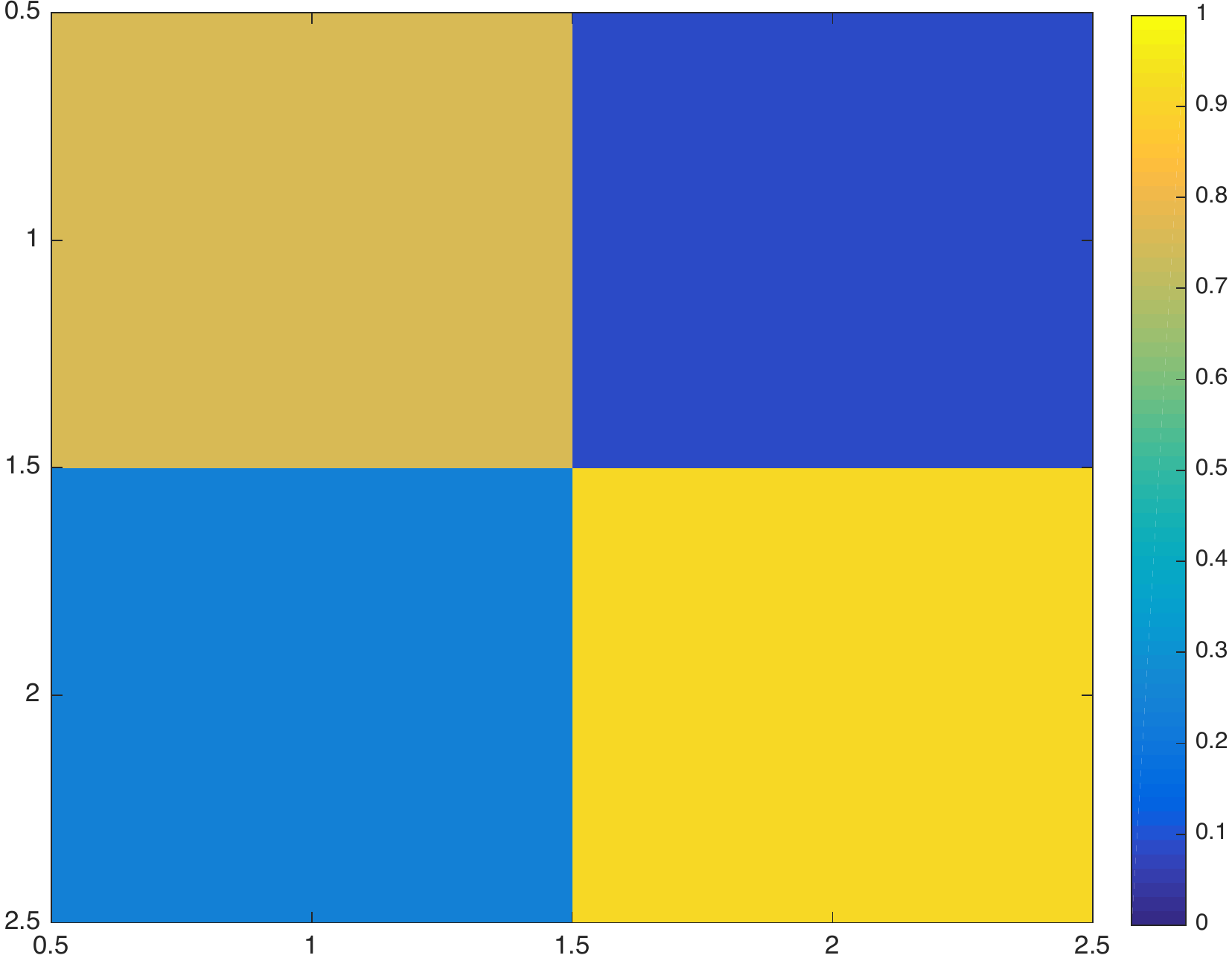}\\[1mm]
  \subcaption{{\GCN}}  
  \end{minipage}
  \begin{minipage}{0.16\textwidth}
  \includegraphics[width=0.99\textwidth,height=0.99\textwidth]{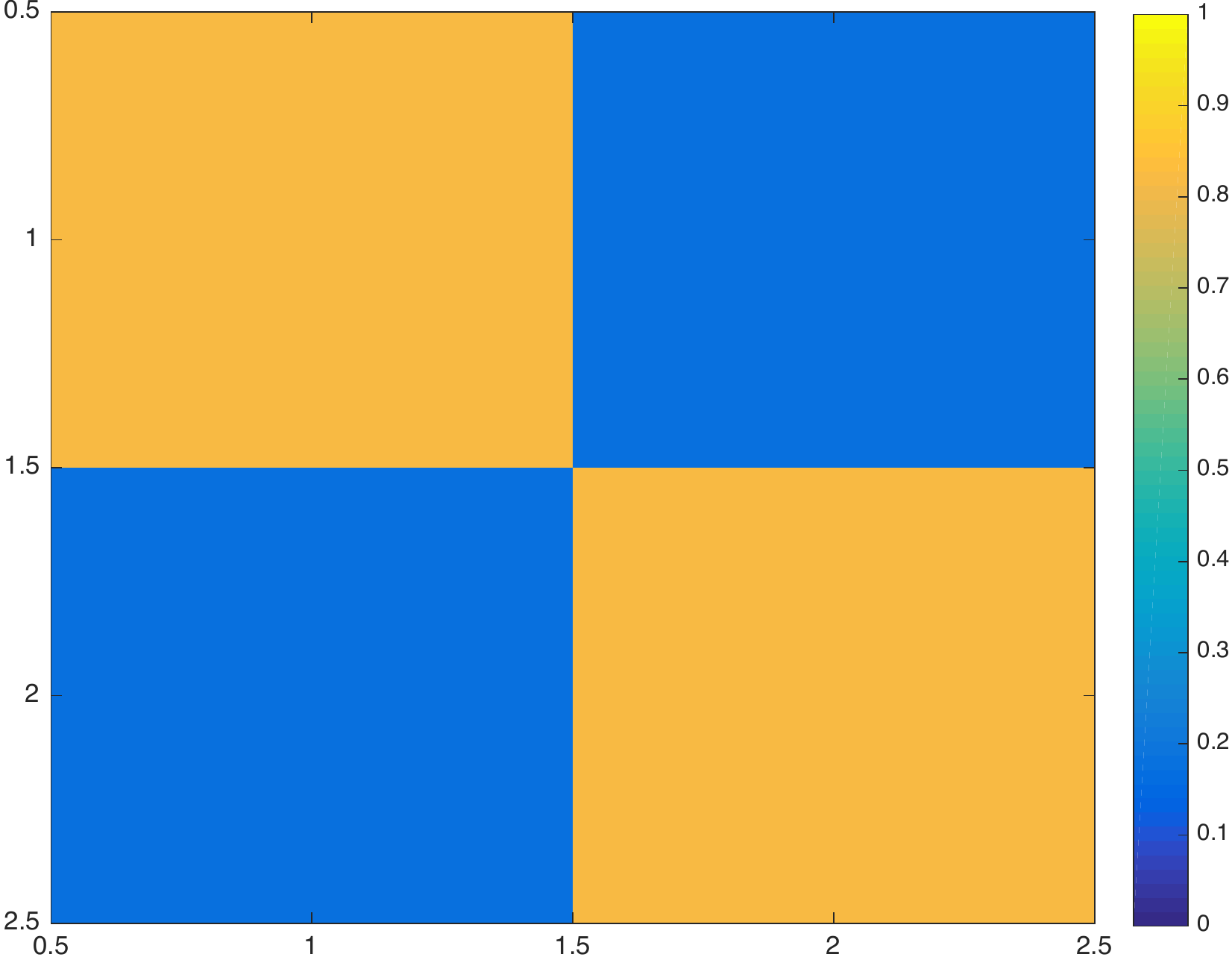}\\[1mm]
  \subcaption{{\LL}} 
  \end{minipage}
  \begin{minipage}{0.16\textwidth}
  \includegraphics[width=0.99\textwidth,height=0.99\textwidth]{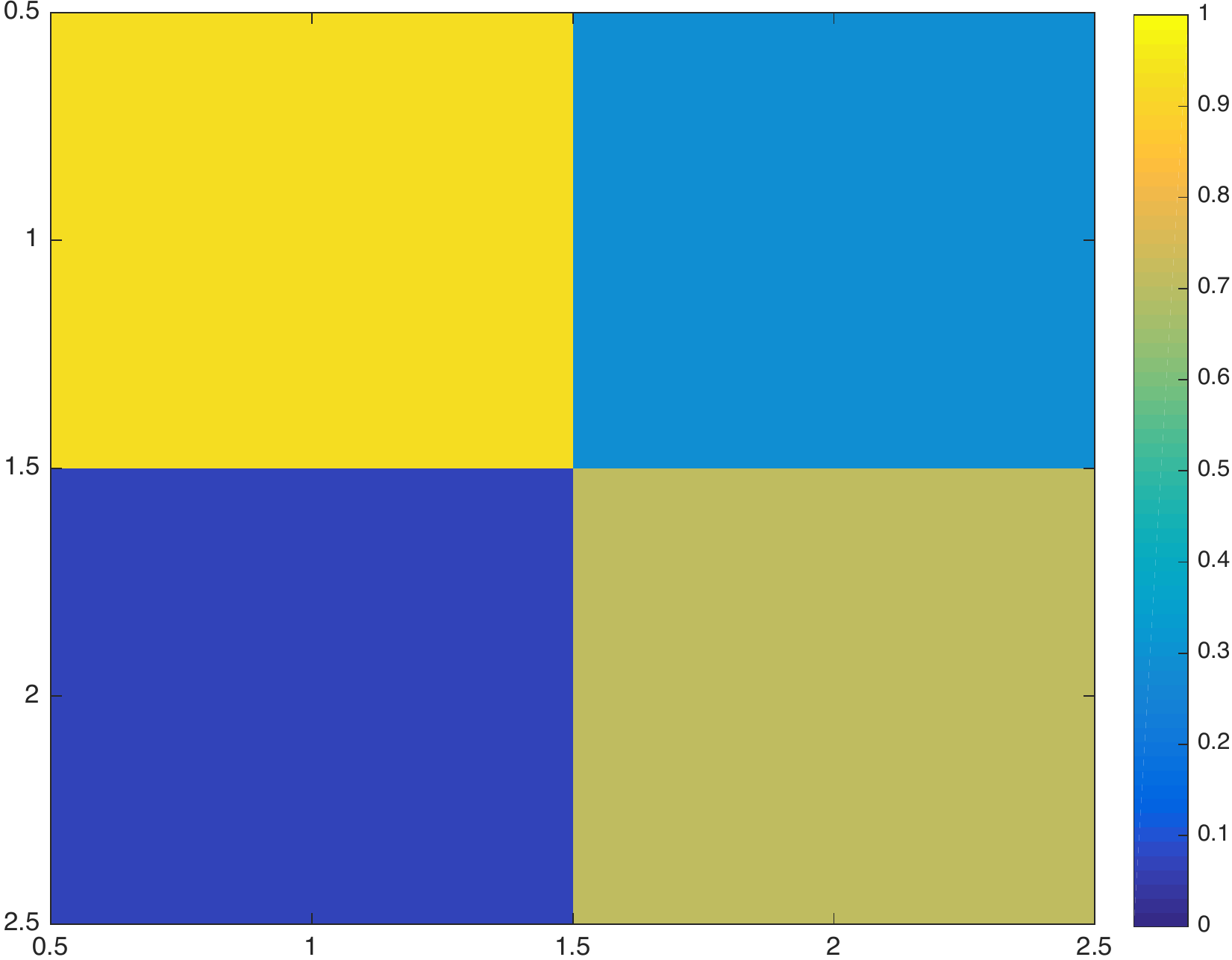}\\[1mm]
  \subcaption{{\Percus}}
  \end{minipage}
\end{minipage}
\caption{Confusion matrices (images) for Proposition dataset. Left 3 (Level 2: 10 classes). Right 3 (Level 1: 2 clases)}\label{fig:CM_Prop}
\end{figure}

In Figures~\ref{fig:CM_CORA} and \ref{fig:CM_Prop},  we present the confusion matrices (as images) for the semi-supervised learning results in Sections~\ref{corasect} and \ref{propdatasect} for datasets CORA and Proposition. We consider 70\% of training data as input and run each clustering of the algorithms ({\GCN}, {\LL}, {\Percus}), which gives different clustering results.  For each method, the confusion matrix $M$ is averaged over 30 runs. From the images of Figures~\ref{fig:CM_CORA} and \ref{fig:CM_Prop}, we can see that for the CORA dataset,  {\GCN} has better classification (visual) results than the other two methods at both level 2 and level 1. For the Proposition dataset,  {\LL} outperforms the other two methods at level 2 while the three mthods performs more or less the same at level 1. Note that since the $F$-measure is computed differently from the class label assignment, the best $F$-measure might not correspond to the best confusion matrix, as measured by misclassification percentage.

We resume the main discussion with a definition of the smoothness classes in terms of the degree of approximation as in \cite{treepap}. Let $0<\rho\le \infty$, $\gamma>0$, and ${\bf a}=\{a_n\}_{n=0}^\infty$ be a sequence of real numbers. We define
\be\label{seqbdef}
\|{\bf a}\|_{\rho,\gamma} := \left\{\begin{array}{ll}
 \disp\left\{\sum_{n=0}^\infty 2^{n\gamma\rho}|a_n|^\rho\right\}^{1/\rho}, & \mbox{if $0<\rho<\infty$,}\\[3ex]
\disp\sup_{n\ge 0} 2^{n\gamma}|a_n|, & \mbox{if $\rho=\infty$.}
\end{array}\right.
\ee
The space of sequences ${\bf a}$ for which $\|{\bf a}\|_{\rho,\gamma}<\infty$ will be denoted by $\seqb_{\rho,\gamma}$. The smoothness class $B_{\rho,\gamma}$ is defined by
\be\label{besovspacedef}
B_{\rho,\gamma}:=\{f\in \mathcal{F} : \{E_{2^n}(f)\}_{n=0}^\infty \in \seqb_{\rho,\gamma}\}.
\ee

In order to develop both function approximation estimates and a wavelet--like characterization of smoothness spaces using the atoms $\tau_j(f)$, we need  an appropriate notion of derivatives. In the case of classical theory of multivariate Fourier series, this is typically done via multipliers; e.g., denoting by $\hat{f}(\k)$ the trigonometric Fourier coefficient, the mixed partial derivative of $f$ is a function $g$ with $\hat{g}(\k)=k_1k_2\hat{f}(k)$, and a spherical partial derivative is a function $g$ with $\hat{g}(\k)=(|\k|^2+1)^{1/2}\hat{f}(\k)$. We have adopted a similar strategy also for tree polynomials in \cite{treepap}.
In the current context, since we don't know the structure of $\Omega$, it is not feasible to define a derivative by means of a fixed multiplier sequence. The following Definition~\ref{derivativedef} gives our substitute for a derivative of order $r$. 

For any subset $S\subseteq \ZZ_+^2$, we define 
$$
\mathcal{E}(S)=\{\k\in\ZZ_+^2 : \m\in S \mbox{ for
some $\m$ with } |\k-\m|_\infty \le 2\}.
$$
 By the restriction of a sequence $h$ to $S$, we mean the sequence whose value at $\k\in S$ is $h(\k)$, and $0$ otherwise. The sequence $h^{[-1]}$ is defined by $h^{[-1]}(\k)=(h(\k))^{-1}$ if $h(\k)\not=0$ and $h^{[-1]}(\k)=0$ otherwise.

\begin{definition}\label{derivativedef}
Let $r\ge 1$ be an integer. A sequence $\mu$ is called a ($\bs g$--)multiplier sequence of order $r$ if $\mu(\k)> 0$ for every $\k\in\mathcal{E}(\Omega)$, and for every integer $j\ge 0$, if $\mu_j$ is the sequence $\mu(\k)$ restricted to $\mathsf{supp}(g_j)$ , then
\be\label{multiplierdef2}
\mathcal{V}(\mu_j)\sim (\mathcal{V}(\mu_j^{[-1]}))^{-1}\sim 2^{jr}
\ee
The derivative of  $f:\GG\to\RR$ (in the sense of $\mu$ and $\bs g$) is a function $\mathcal{D}(f): \GG\to\RR$ such that
$\widehat{\mathcal{D}(f)}(\k)=\mu(\k)\hat{f}(\k)$, $\k\in\ZZ_+^2$, if such a function exists.
\end{definition}

In the rest of this paper, we will fix $r$ and $\mu$. 

The $K$--functional in the theory of function approximation on $\GG$ is now defined for $f\in \mathcal{F}$ by
\be\label{kfuncdef}
K(f,\delta)=\inf\{\|f-g\|+\delta^r\|\mathcal{D}g\| : g, \mathcal{D}g\in\mathcal{F}\}.
\ee

The analogue of \cite[Theorem~5.2, Theorem~5.3]{treepap} is the following.

\begin{theorem}\label{tenswavetheo}
Let $0<\rho\le \infty$ and $\gamma>0$. Let $\mathcal{V}(H_n)\lesssim 1$ for all $n\ge 1$, and $f\in \mathcal{F}$. The following statements are equivalent.
\begin{enumerate}
\item $f\in B_{\rho,\gamma}$.
\item $\{\|\tau_j(f)\|\}_{j=0}^\infty \in \seqb_{\rho,\gamma}$.
\item $\{\|f-\sigma_n(f)\|\}_{j=0}^\infty \in \seqb_{\rho,\gamma}$.
\item $\{K(f,2^{-j})\}\in \seqb_{\rho,\gamma}$.
\end{enumerate}
\end{theorem}

\bhag{Proof of the theorems in Section~\ref{diharmsect}.}\label{pfsect}
In the sequel, we find it convenient to extend the sequences on $\Omega$ to sequences on $\ZZ_+^2$ by setting their values to be $0$ outside $\Omega$.  This will facilitate the use of summation indices and 
relationships such as \eref{leibnitz} and \eref{leibineq} provided we keep in mind that the validity of \eref{partunity} is assumed only on $\Omega$. Similarly, functions on $\GG$ are assumed to be extended to $I^2$ by setting them equal to $0$ outside of $\GG$.

\noindent
\textsc{Proof of Theorem~\ref{tensfourtheo}.}
In this proof, let 
$$
s_{(n_1,n_2)}(f)(x_1,x_2)=\sum_{k_1=0}^{n_1}\sum_{k_2=0}^{n_2}\hat{f}(k_1,k_2)\psi_{k_1,1}(x_1)\psi_{k_2,2}(x_2), \qquad n_1, n_2=0,1,\cdots,
$$
where $\hat{f}(k_1,k_2)=0$ if $(k_1,k_2)\not\in\Omega$. Then \cite[Theorem~5.1]{treepap} regarding the uniform boundedness of the Fourier partial sums in the univariate case leads to
\be\label{pf1eqn1}
\|s_\k(f)\|\lesssim 1, \qquad \k\in\ZZ_+^2.
\ee 
Next, we observe that $\hat{f}(\k)\psi_\k(\x)=\bs\Delta s_\k(f)(\x)$ for all $\k\in\ZZ_+^2$. A summation by parts shows that
$$
\sum_{\k\in\ZZ_+^2} h(\k)\hat{f}(\k)\psi_\k(\x)=\sum_\k \bs\Delta h(k_1-1,k_2-1)s_\k(f)(\x),
$$
where $h(k_1,k_2)=0$ if either $k_1$ or $k_2$ is negative. The estimate \eref{tensgenseqfourbd} follows from \eref{pf1eqn1} and the definition of $\mathcal{V}(h)$. \qed

\noindent
\textsc{Proof of Theorem~\ref{tensgoodapproxtheo}.} The first estimate in \eref{tensgoodapprox} is obvious from the definition. The second estimate is also obvious in light of Theorem~\ref{tensfourtheo} and the definition if $n< m^*$.  In the remainder of this proof, let $n\ge m^*$, and
$$
\Omega_n=\{\k : H_n(\k)\not=0\}, \qquad n=0,1,\cdots.
$$
If $\k\in\Omega_{n-m^*}$, then there exists $j\le n-m^*$ such that $g_j(\k)\not=0$. Since $g_j(\k)g_{j'}(\k)=0$ if $|j-j'|>m^*$, this implies that $g_{j'}(\k)=0$ for all $j'>n$. Consequently, \eref{partunity} shows that $H_n(\k)=\sum_{j=0}^n g_j(\k)=\sum_{j\in\Omega} g_j(\k)=1$. 
So, if $P=\sum_{\k\in\Omega_{n-m^*}} \hat{P}(\k)\psi_\k \in \mathbb{P}_{n-m^*}$, then for all $\x\in\GG$,
$$
\sigma_n(P)(\x)=\sum_{\k}H_n(\k)\hat{P}(\k)\psi_\k(\x)=\sum_{\k\in\Omega_{n-m^*}}H_n(\k)\hat{P}(\k)\psi_\k(\x)= \sum_{\k\in\Omega_{n-m^*}}\hat{P}(\k)\psi_\k(\x)=P(\x).
$$
Therefore, Theorem~\ref{tensfourtheo} leads to 
$$
\|f-\sigma_n(f)\|\le \|f-P\|+\|\sigma_n(f-P)\|\lesssim\|f-P\|,
$$
and hence, to the second estimate of \eref{tensgoodapprox}.

The second estimate in \eref{tensgoodapprox} can be rewritten in the form
$$
\left|f-\sum_{j=0}^n\tau_j(f)\right|=\|f-\sigma_n(f)\|\lesssim E_{n-m^*}(f).
$$
Since $f\in\mathcal{F}$, $E_{n-m^*}(f)\to 0$ as $n\to\infty$, this is equivalent to \eref{tenslpseries} in the sense of uniform convergence. This proves part (b). 

To prove part (c), we observe that since $\nu^*$ is a probability measure, the system $\{\psi_\k\}$ is a fundamental system for the $L^2(\GG,\nu^*)$ closure of $\mathcal{F}$. Therefore, for $f\in\mathcal{F}$, Parseval identity holds, and we obtain
\be\label{pf1eqn2}
\int_\GG |f(\x)|^2d\nu^*(\x)=\sum_{\k\in\ZZ_+^2}|\hat{f}(\k)|^2 = \sum_{j=0}^\infty \sum_{\k\in\ZZ_+^2}g_j(\k)|\hat{f}(\k)|^2.
\ee
Since $g_j(\k)^2\le g_j(\k)$ for all $j$ and $\k$, this shows using Parseval identity again that
\be\label{pf1eqn3}
\sum_{j=0}^\infty \int_\GG |\tau_j(f)(\x)|^2d\nu^*(\x)=\sum_{j=0}^\infty \sum_{\k\in\ZZ_+^2}g_j(\k)^2|\hat{f}(\k)|^2\le \int_\GG |f(\x)|^2d\nu^*(\x).
\ee

Since $g_j(\k)g_m(\k)=0$ if $|j-m|>m^*$, it is easy to verify using the definitions that
$$
\int_\GG \tau_j(f)(\x)\tau_m(f)(\x)d\nu^*(\x)=0, \qquad |j-m|>m^*.
$$
Using \eref{tenslpseries}, we see that
$$
\int_\GG |f(\x)|^2d\nu^*(\x)=\sum_{j=0}^\infty\sum_{m=0}^\infty \int_\GG \tau_j(f)(\x)\tau_m(f)(\x)d\nu^*(\x)= \sum_{j=0}^\infty\sum_{m=\max(j-m^*,0)}^{j+m^*} \int_\GG \tau_j(f)(\x)\tau_m(f)(\x)d\nu^*(\x).
$$
An application of Schwarz inequality and Parseval theorem then lead to
$$
 \int_\GG |f(\x)|^2d\nu^*(\x)\lesssim \sum_{j=0}^\infty \int_\GG |\tau_j(f)(\x)|^2d\nu^*(\x).
 $$
Together with \eref{pf1eqn3}, this completes the proof of \eref{frameprop}. \qed

\begin{theorem}\label{fundaineqtheo}
{\rm (a)} For $n\ge 0$ and $P\in \mathbb{P}_n$,
\be\label{bernineq}
\|\mathcal{D} P\| \lesssim2^{nr}\|P\|.
\ee
{\rm (b)} If $f\in\mathcal{F}$ and $\mathcal{D}^r f\in \mathcal{F}$, then for $n\ge 0$,
\be\label{favardineq}
E_n(f)\lesssim 2^{-nr}\|\mathcal{D} f\|.
\ee
\end{theorem} 

\begin{Proof}\ 
We observe first that for $j=0,1, \cdots$ and $\k\in\ZZ_+^2$, $\widehat{\tau_j(f)}(\k)=g_j(\k)\hat{f}(\k)$. Hence, Theorem~\ref{tensfourtheo} and the conditions \eref{multiplierdef2} together imply that
\be\label{pf2eqn1}
\|\mathcal{D}\tau_j(f)\|=\left\|\sum_\k \mu(\k)g_j(\k)\hat{f}(k)\psi_\k\right\| =\left\|\sum_\k \mu_j(\k)g_j(\k)\hat{f}(k)\psi_\k\right\|\lesssim 2^{jr}\|\tau_j(f)\|\lesssim 2^{jr}\|f\|,
\ee
and similarly,
\be\label{pf2eqn2}
\|\tau_j(f)\|=\left\|\sum_\k \mu_j^{[-1]}(\k)\mu_j(\k)g_j(\k)\hat{f}(\k)\psi_\k\right\| \lesssim 2^{-jr}\left\|\sum_\k g_j(\k)\mu(\k)\hat{f}(\k)\psi_\k\right\|=2^{-jr}\|\tau_j(\mathcal{D}f)\|\lesssim 2^{-jr}\|\mathcal{D}f\|.
\ee

To prove part (a), we observe  in view of \eref{pf2eqn1} that
$$
\|\mathcal{D}P\|=\|\mathcal{D}\sigma_{n+m^*}(P)\|=\left\|\sum_{j=0}^{n+m^*} \mathcal{D}\tau_j(P)\right| \lesssim \sum_{j=0}^{n+m^*} 2^{jr}\|P\|\lesssim 2^{nr}\|P\|.
$$
This proves \eref{bernineq}.

To prove part (b), we use \eref{tenslpseries} and \eref{pf2eqn2} to deduce that
$$
E_n(f)\le \left\|\sum_{j=n+1}^\infty \tau_j(f)\right\|\le \sum_{j=n+1}^\infty \|\tau_j(f)\|\lesssim \|\mathcal{D}f\|\sum_{j=n+1}^\infty 2^{-jr} \lesssim 2^{-nr}\|\mathcal{D}f\|.
$$
This proves \eref{favardineq}.
\end{Proof}

\noindent
\textsc{Proof of Theorem~\ref{tenswavetheo}.} Using Theorem~\ref{fundaineqtheo}, the proof of  Theorem~\ref{tenswavetheo} follows standard arguments. For the equivalence of the items 1, 2, 3, these arguments are exactly the same as those in \cite[Theorem~4]{fasttour} (with different notation). The equivalence of items 1 and 4 is shown using the arguments in \cite[Theorem~9.1 in Section~7.9, also Chapter 6.7]{devlorbk}. We omit the details. \qed

\bhag{Acknowledgments}
We thank  Professors Percus and Hunter at Claremont Graduate University and Claremont McKenna College respectively for many useful discussions as well as their help in securing the Proposition data set, which was sent to us by Dr. Linhong Zhu at USC Information Sciences Institute in Marina Del Ray, California. We thank Dr. Garcia--Cardona for giving us a C code for the algorithm {\Percus}. 


\end{document}